\documentclass[12pt,a4paper]{article}
\usepackage{amsmath,amsfonts,amssymb,amscd}

\def\al{\alpha}

\def\gr{\operatorname{gr}}
\def\Ker{\operatorname{Ker}}
\def\Ad{\operatorname{Ad}}
\def\ad{\operatorname{ad}}

\def\id{\operatorname{id}}
\def\ev{\operatorname{ev}}
\def\red{\operatorname{red}}
\def\d{\operatorname{d}}
\def\GL{\operatorname{GL}}

\def\pr{\operatorname{pr}}

\def\pt{\operatorname{pt}}
\def\Hom{\operatorname{Hom}}

\def\Aut{\operatorname{Aut}}
\def\Lie{\operatorname{Lie}}

\newcounter{th}
\def\t{\refstepcounter{th}{\bf \noindent{Theorem} \arabic{th}. }}

\newcounter{le}
\def\l{\refstepcounter{le}{\bf \noindent{Proposition} \arabic{le}. }}

\newcounter{lem}
\def\lem{\refstepcounter{lem}{\bf \noindent{Lemma} \arabic{lem}. }}

\newcounter{de}
\def\de{\refstepcounter{de}{\bf \noindent{Definition} \arabic{de}. }}

\begin{document}

\begin{center}
\large{\bf ON COMPLEX LIE SUPERGROUPS AND SPLIT HOMOGENEOUS SUPERMANIFOLDS
\footnote{Work supported by SFB $\mid$ TR 12 and by the Russian Foundation
for Basic Research (grant no. 07-01-00230).
 } }
\end{center}

\begin{center}
     E.G. Vishnyakova
\end{center}

\noindent\textsc{Abstract.} It is well known that the category of
real Lie supergroups is equivalent to the category of the so called
(real) Harish-Chandra pairs, see \cite{Bern, Kostant, kosz}. That
means that a Lie supergroup depends only on the underlying Lie group
and its Lie superalgebra with certain compatibility conditions. More
precisely, the structure sheaf of a Lie supergroup and the
supergroup morphisms can be explicitly described in terms of the
corresponding Lie superalgebra. In this paper, we give a proof of
this result in the complex-analytic case. Furthermore, if
$(G,\mathcal{O}_G)$ is a complex Lie supergroup and $H\subset G$ is
a closed Lie subgroup, i.e. it is a Lie subsupergroup of
$(G,\mathcal{O}_G)$ and its odd dimension is zero, we show that the
corresponding homogeneous supermanifold $(G/H,\mathcal{O}_{G/H})$ is
split. In particular, any complex Lie supergroup is a split
supermanifold.

It is well known that a complex homogeneous supermanifold may be non-split
(see, e.g., \cite{oniosp}). We find here necessary and sufficient conditions
for a complex homogeneous supermanifold to be split.

\bigskip

\begin{center}
{\bf 1. Preliminaries}
\end{center}

We will use the word "supermanifold" \,in the sense of Berezin --
Leites (see \cite{BL, ley}).  All the time, we will be interested in
the real or complex-analytic version of the theory, denoting by
$\Bbb K$ the ground field $\Bbb R$ or $\Bbb C$. Let
$(M,\mathcal{O}_M)$ be a supermanifold. The underlying complex
manifold $M$ is called the {\it reduction} of $(M,\mathcal{O}_M)$.
We denote by $\mathcal{J}_M\subset\mathcal{O}_M$ the subsheaf of
ideals generated by odd elements of the structure sheaf. The sheaf
$\mathcal{O}_M/\mathcal{J}_M$ is naturally identified with the
structure sheaf $\mathcal{F}_M$ of $M$. The natural homomorphism
$\mathcal{O}_M\to\mathcal{F}_M$ will be denoted by $f\mapsto
f_{\red}$.
 A morphism
$\phi:(M,\mathcal{O}_M)\to (N,\mathcal{O}_N)$ of supermanifolds will
be denoted by $\phi = (\phi_{\red},\phi^*)$, where $\phi_{\red}:
M\to N$ is the corresponding mapping of the reductions and
$\phi^*:\mathcal{O}_N\to(\phi_{\red})_*(\mathcal{O}_M)$ is the
homomorphism of the structure sheaves. We denote by
$\mathfrak{v}(M,\mathcal{O}_M)$ the Lie superalgebra of vector
fields on $(M,\mathcal{O}_M)$. If $x\in M$ and $\mathfrak m_x$ is
the maximal ideal of the local superalgebra $(\mathcal O_M)_x$, then
the vector superspace $T_x(M,\mathcal O_M)=(\mathfrak m_x/\mathfrak
m_x^2)^*$ is the tangent space to $(M,\mathcal O_M)$ at $x\in M$.
From the inclusions $v(\mathfrak m_x)\subset (\mathcal O_M)_x$ and
$v(\mathfrak m_x^2)\subset \mathfrak m_x$, where $v\in \mathfrak
v(M,\mathcal O_M)$, it follows that $v$ induces an even linear
mapping $\ev_x(v):\mathfrak m_x/\mathfrak m_x^2\to (\mathcal
O_M)_x/\mathfrak m_x\simeq \Bbb C$. In other words, $\ev_x(v)\in
T_x(M,\mathcal O_M)$, and so we obtain an even linear map
\begin{equation}\label{ev}
\ev_x:\mathfrak v(M,\mathcal O_M)\to T_x(M,\mathcal O_M).
\end{equation}

Let us take $Y_x\in T_x(M,\mathcal O_M)$. There is a neighborhood
$(U,\mathcal O_M)$ of the point $x$ and a vector field $Y\in
\mathfrak{v}(U,\mathcal O_M)$ such that $\ev_x(Y)=Y_x$. We may
regard $Y_x$ as a linear function on $(\mathcal O_M)_x$. Namely,
 $Y_x(f_x):=(Y(f_x))_{\red}(x)$, where
$f_x\in (\mathcal O_M)_x$. It is easy to verify that this definition
doesn't depend on the choice of $Y$.

Let $(M,\mathcal{F}_M)$ be a complex-analytic, smooth or
real-analytic manifold and let $\mathbf{E}_M$ be a (holomorphic,
smooth or real-analytic) vector bundle over $M$. Denote by
$\mathcal{E}_M$ the sheaf of (holomorphic, smooth or real-analytic)
sections of $\mathbf{E}_M$. We get the supermanifold
$(M,\bigwedge_{\mathcal{F}_M}\mathcal{E}_M)$ of the corresponding
class. A supermanifold $(M,\mathcal{O}_M)$ is called split if
$\mathcal{O}_M\simeq \bigwedge_{\mathcal{F}_M}\mathcal{E}_M$ for a
certain vector bundle $\mathbf{E}_M$ of the corresponding class. It
is known that any real (smooth or real analytic) supermanifold is
split.

We may consider the supermanifold $(\pt,\mathbb{K})$ of dimension
$(0|0)$, where $\pt$ is a point. If $(M,\mathcal{O}_M)$ is an
arbitrary supermanifold then for any point $x\in M$ we denote by
$\delta_x:(\pt,\mathbb{K})\to (M,\mathcal{O}_M)$ the morphism,
defined in the following way:
$$
\begin{array}{l}
(\delta_x)_{\red}(\pt)=x,\,\,\,\,
 \delta_x^*(f) = \left\{
  \begin{array}{ll}
    f_{\red}(x), & \hbox{if $x\in U$;} \\
    0, & \hbox{if $x\notin U$,}
  \end{array}
\right.
\end{array}
$$
where $f\in \mathcal{O}_M(U)$ and $U\subset M$ is open.

\medskip

\de\label{def Liesgroup} A {\it Lie supergroup} is a group object in
the category of supermanifolds, i.e., a supermanifold $(G,\mathcal
O_G)$, for which the following three morphisms are defined:
$\mu:(G,\mathcal O_G)\times (G,\mathcal O_G)\to (G,\mathcal O_G)$
(the multiplication morphism), $\iota:(G,\mathcal O_G)\to
(G,\mathcal O_G)$ (the inversion morphism),
$\varepsilon:(\pt,\mathbb K)\to (G,\mathcal O_G)$ (the identity
morphism). Moreover, these morphisms should satisfy the usual
conditions, modeling the group axioms:

\begin{enumerate}
  \item $\mu\circ (\mu\times \id) = \mu\circ (\id\times \mu)$;
  \item $\mu\circ (\varepsilon\times \id)=\id$, $\mu\circ (\id\times
  \varepsilon)=\id$;
  \item $\mu\circ (\id\times \iota)\circ \operatorname{diag}=\varepsilon$, $\mu\circ (\iota\times
  \id)\circ \operatorname{diag}=\varepsilon$, where $\operatorname{diag}:
  (G,\mathcal O_G)\to (G,\mathcal O_G)\times (G,\mathcal O_G)$ is the diagonal morphism.
\end{enumerate}
The underlying manifold $G$ of a Lie supergroup is a (real or
complex) Lie group. The element $e=\varepsilon_{\red}(\pt)$ is the
identity element of $G$. Let $(G,\mathcal{O}_G)$,
$(H,\mathcal{O}_H)$ be two Lie supergroups and $\mu_G$, $\mu_H$ the
respective multiplication morphisms. A morphism
$\Psi:(G,\mathcal{O}_G) \to (H,\mathcal{O}_H)$ is called a {\it
homomorphism of Lie supergroups} if $\Psi\circ \mu_G=\mu_H \circ
(\Psi\times \Psi)$. The corresponding mapping $\Psi_{\red}:G\to H$
is a homomorphism of Lie groups.

\medskip

\de\label{def sgroup action} An {\it action of a Lie supergroup
$(G,\mathcal O_G)$ on a supermanifold} $(M,\mathcal O_M)$ is a
morphism $\nu:(G,\mathcal O_G)\times (M,\mathcal O_M)\to (M,\mathcal
O_M)$, such that the following conditions hold:
\begin{itemize}
  \item $\nu \circ (\mu\times \id)=\nu\circ (\id\times \nu)$;
  \item $\nu\circ (\varepsilon\times \id)=\id$.
\end{itemize}
In this case $\nu_{\red}$ is the action of $G$ on $M$.

We will denote by $\mathfrak{g}$ the Lie superalgebra of
$(G,\mathcal O_G)$. By definition $\mathfrak{g}$ is the subalgebra
of $\mathfrak{v}(G,\mathcal{O}_G)$, consisting of all right
invariant vector fields on $(G,\mathcal{O}_G)$. (A vector field $Y$
on $(G,\mathcal{O}_G)$ is called right invariant if $(Y\otimes
\id)\circ \mu^*= \mu^* \circ Y$.) It is well known that any right
invariant vector field $Y$ has the form
\begin{equation}
\label{left inv vect field}
 Y=(X\otimes \id)\circ \mu^*
\end{equation}
 for a certain $X\in T_e(G,\mathcal{O}_G)$ and the map $X\mapsto (X\otimes \id)\circ \mu^*$ is an
isomorphism of the vector space $T_e(G,\mathcal{O}_G)$ onto
$\mathfrak{g}$, see \cite{Var}, Theorem $7.1.1$. We will identify
$\mathfrak{g}$ and $T_e(G,\mathcal{O}_G)$ using this isomorphism.

Let $\nu=(\nu_{\red},\nu^*):(G,\mathcal O_G)\times (M,\mathcal
O_M)\to (M,\mathcal O_M)$ be an action. Then there is a homomorphism
of the Lie superalgebras $\overline{\nu}:\mathfrak g\to \mathfrak
v(M,\mathcal O_M)$, given by the formula
\begin{equation}\label{hom}
 X\mapsto
(X\otimes\, \id)\circ \nu^*.
\end{equation}

As in \cite{onipi}, we use the following definition of a transitive
action.

\medskip

\de\label{def homogen supermnf}   An action $\nu$ is called {\it
transitive} if $\nu_{\red}$ is transitive and the mapping $\ev_{x}
\circ \overline{\nu}$ is surjective for all $x\in M$. (The map
$\ev_{x}$ is given by (\ref{ev}).) In this case the supermanifold
$(M,\mathcal O_M)$ is called $(G,\mathcal O_G)$-{\it homogeneous}.
 A supermanifold $(M,\mathcal O_M)$
is called {\it homogeneous}, if it possesses a transitive action of
a certain Lie supergroup.

Let us consider the following compositions of the morphisms for any
$g\in G$:
$$
\begin{array}{l} l_g:(G,\mathcal{O}_G)=(g,\mathbb{K})\times
(G,\mathcal{O}_G)\stackrel{\delta_g\times \id}{\longrightarrow}
(G,\mathcal{O}_G)\times
(G,\mathcal{O}_G)\stackrel{\mu}{\to}(G,\mathcal{O}_G),\\
r_g:(G,\mathcal{O}_G)=(G,\mathcal{O}_G)\times
(g,\mathbb{K})\stackrel{\id\times \delta_g}{\longrightarrow}
(G,\mathcal{O}_G)\times
(G,\mathcal{O}_G)\stackrel{\mu}{\to}(G,\mathcal{O}_G).
\end{array}
$$
They are called the left and the right translation by $g$
respectively. Denote $\omega_g:=l_g\circ r_{g^{-1}}$, $g\in G$. The
formula $\Ad_G(g):= (\d \omega_g)_e$ defines a representation
$\Ad_G:G\to \Aut(\mathfrak{g})$, called the {\it adjoint
representation of the Lie group $G$} in $\mathfrak{g}$.

Let $(M,\mathcal{O}_M)$ be a supermanifold. A {\it subsupermanifold}
of $(M,\mathcal{O}_M)$ is a supermanifold $(N,\mathcal{O}_N)$
together with a morphism $\varphi: (N,\mathcal{O}_N)\to
(M,\mathcal{O}_M)$ such that $\varphi_{\red}:N\to M$ is a
homeomorphism on the subset $\varphi_{\red}(N)\subset M$ endowed
with the induced topology and $(\d\varphi)_p$ is injective at every
point $p\in M$. In this case we will sometimes use the notation
$(M,\mathcal O_{M}) \subset (N,\mathcal O_{N})$.

Let $(G,\mathcal O_G)$ be a Lie supergroup. We say that a
subsupermanifold $\varphi: (H,\mathcal O_H)\to (G,\mathcal O_G)$ is
a {\it Lie subsupergroup} in $(G,\mathcal O_G)$ if $(H,\mathcal
O_H)$ possesses a Lie supergroup structure, such that $\varphi$ is a
homomorphism of the Lie supergroups. In this case we identify the
Lie superalgebra $\mathfrak{h}$ of $(H,\mathcal O_H)$ with the Lie
subsuperalgebra $(\d\varphi)_e(\mathfrak{h})\subset\mathfrak{g}$.

Let us introduce the category of (super) Harish-Chandra pairs (see
\cite{Bern}). A {\it Harish-Chandra pair} is a pair
$(G,\mathfrak{g})$ that consists of a Lie group $G$ and a Lie
superalgebra $\mathfrak{g}=\mathfrak{g}_{\bar
0}\oplus\mathfrak{g}_{\bar 1}$, where $\mathfrak{g}_{\bar 0}$ is the
Lie algebra of $G$, provided with a representation $\al_G$ of $G$ in
$\mathfrak{g}$ such that
\begin{itemize}
  \item $\al_G$ preserves the parity and induces the adjoint representation
of $G$ in $\mathfrak{g}_{\bar 0}$,
  \item the differential $(\d\al_G)_e$ at the identity $e\in G$ coincides with
the adjoint representation $\ad$ of $\frak g_{\bar 0}$ in $\frak g$.
\end{itemize}

Let $(G,\mathfrak{g})$ and $(H,\mathfrak{h})$ be two Harish-Chandra
pairs. A {\it morphism} of $(G,\mathfrak{g})$ to $(H,\mathfrak{h})$
is a pair of homomorphisms $\Phi:G\to H$, $\varphi: \mathfrak{g}\to
\mathfrak{h}$ with the following compatibility conditions:
\begin{itemize}
  \item $(\d\Phi)_e = \varphi\mid_{\mathfrak{g}_{\bar 0}}$,
  \item $\varphi\circ\al_G(g)=\al_H(\Phi(g))\circ \varphi$ for all
  $ g\in G$.
\end{itemize}

It is clear how to associate a Harish-Chandra pair to a given Lie
supergroup $(G,\mathcal{O}_G)$. Indeed, we may take the underlying
Lie group $G$ with the Lie superalgebra $\mathfrak{g}$ of
$(G,\mathcal{O}_G)$ equipped with the adjoint representation $\al_G
= \Ad_G$. Furthermore, if $\Psi:(G,\mathcal{O}_G) \to
(H,\mathcal{O}_H)$ is a homomorphism of Lie supergroups, then
$(\Psi_{\red}, (\d \Psi)_e)$ is a morphism of the Harish-Chandra
pairs $(G,\mathfrak{g})\to(H,\mathfrak{h})$. This correspondence is
a functor from the category of Lie supergroups to the category of
Harish-Chandra pairs. From Theorem $3.5$ and Remark 3.5.2 in
\cite{Kostant} it can be deduced that this functor is an equivalence
of categories in the real case. The proof in \cite{Kostant} uses the
fact that $C^{\infty}$-supermanifold can be reconstructed from the
algebra of global sections of its structure sheaf (see \cite[Remark
2.14.2]{Kostant}). Since such reconstruction is in general
impossible for holomorphic supermanifolds, this argument doesn't
seem to immediately carry over to the holomorphic case. We will give
a different proof of the equivalence that works both in the real and
holomorphic cases.

Let us denote the category of Harish-Chandra pairs by $\verb"HCP"$
and the category of Lie supergroups by $\verb"SLG"$.

\bigskip


\begin{center}
  {\bf 2.  Equivalence between $\verb"HCP"$ and $\verb"SLG"$}
\end{center}

In this section we will prove that the categories $\verb"HCP"$ and
$\verb"SLG"$ are equivalent. We denote by $\operatorname{Ob}\,C$ the
set of objects of a category $C$ and by $\Hom(X,Y)$ the set of
morphisms $X\to Y$ for two objects $X,Y\in \operatorname{Ob}\,C$.
First, we shall describe a functor $F$ from the category
$\verb"HCP"$ to $\verb"SLG"$ that was constructed by Koszul in
\cite{kosz}. Further, we show that for any object $Y\in
\operatorname{Ob}\,\verb"SLG"$ there exists $X\in
\operatorname{Ob}\,\verb"HCP"$ such that $F(X)$ is isomorphic to
$Y$. Finally, we prove that $F: \Hom(X,Y)\to\Hom(F(X), F(Y))$ is a
bijection for every $X,Y\in \operatorname{Ob}\,\verb"HCP"$. This
will imply that $F$ determines an equivalence of our categories (see
\cite{TSH}).

\medskip

\noindent{\it 2.1 The construction of $F$.} If a (real or complex)
Harish-Chandra pair $(G,\mathfrak{g})$ is given, then we can
construct a Lie supergroup in the following way (see \cite{BagSta,
kosz}). Let $\mathfrak{U}(\mathfrak{g})$ be the universal enveloping
superalgebra of $\mathfrak{g}$ (see \cite{Scheu}). It is clear that
$\mathfrak{U}(\mathfrak{g})$ is a $\mathfrak{U}(\mathfrak{g}_{\bar
0})$-module, where $\mathfrak{U}(\mathfrak{g}_{\bar 0})$ is the
universal enveloping algebra of $\mathfrak{g}_{\bar 0}$. The natural
action of $\mathfrak{g}_{\bar 0}$ on the sheaf $\mathcal{F}_G$ gives
rise to a structure of $\mathfrak{U}(\mathfrak{g}_{\bar 0})$-module
on $\mathcal{F}_G(U)$ for any open set $U\subset G$. Putting
$$
\widehat{\mathcal{O}}_G(U) =
\Hom_{\mathfrak{U}(\mathfrak{g}_{\bar 0})}(\mathfrak{U}(\mathfrak{g}),
\mathcal{F}_G(U))
$$
for every open $U\subset G$, we get a sheaf $\widehat{\mathcal{O}}_G$ of
$\mathbb{Z}_2$-graded vector spaces (here we assume that the functions from
$\mathcal{F}_G(U)$ are even).

As a consequence of the graded version of Theorem of
Poincar\'{e}-Birkhoff-Witt, we obtain that $
\mathfrak{U}(\mathfrak{g}_{\bar 0})\otimes
\bigwedge(\mathfrak{g}_{\bar 1})\simeq \mathfrak{U}(\mathfrak{g})$
as $\mathfrak{U}(\mathfrak{g}_{\bar 0})$-modules (see \cite{kosz,
Scheu}). The isomorphism is given by the formula $X\otimes Y\mapsto
X \cdot \gamma(Y) $, where
\begin{equation}
\label{isomorphism}
 \gamma: \bigwedge(\mathfrak{g}_{\bar 1})\to
\mathfrak{U}(\mathfrak{g}),\quad X_1\wedge \cdots \wedge X_r\mapsto
\frac{1}{r!}\sum_{\sigma\in S_r}(-1)^{|\sigma|} X_{\sigma(1)}\cdots
X_{\sigma(r)}.
\end{equation}
The enveloping superalgebra $\mathfrak{U}(\mathfrak{g})$ has a Hopf
superalgebra structure (see \cite{Scheu}). Indeed, the map
$$
\mathfrak{g} \to \mathfrak{U}(\mathfrak{g})\otimes
\mathfrak{U}(\mathfrak{g}), \,\,\,X\mapsto X\otimes 1 + 1\otimes X;\\
$$
can be extended to a comultiplication map $\bigtriangleup:
\mathfrak{U}(\mathfrak{g})\to \mathfrak{U}(\mathfrak{g})\otimes
\mathfrak{U}(\mathfrak{g})$, and the antipode map
$S:\mathfrak{U}(\mathfrak{g})\to \mathfrak{U}(\mathfrak{g})$ is
given by
$$
S(X)= -X,\,\,\,\, S(1)= 1,\,\,\,\,S(Y\cdot
Z)=(-1)^{p(Y)p(Z)}S(Z)\cdot S(Y),
$$
where $X\in \mathfrak{g}$, $Y,\, Z\in \mathfrak{U}(\mathfrak{g})$
and $p(V)$ is the parity of $V$. We can define a multiplication in
each $\widehat{\mathcal{O}}_G(U)$, where $U\subset G$ is open, by
$$
f_1\cdot f_2:=\operatorname{Mult}_{\mathcal{F}_G} \circ (f_1\otimes
f_2)\circ \bigtriangleup.
$$
Here $f_1,f_2\in \widehat{\mathcal{O}}_G(U)$ and by
$\operatorname{Mult}_{\mathcal{F}_G}$ is denoted the product in the
sheaf $\mathcal{F}_G$. Note that for homogeneous $X,Y\in
\mathfrak{U}(\mathfrak{g})$ and $f_1,f_2\in
\widehat{\mathcal{O}}_G(U)$ we have
\begin{equation}\label{tensorprodukt}
(f_1\otimes f_2)(X\otimes Y) = (-1)^{p(f_2)p(X)} f_1(X)\otimes
f_2(Y).
\end{equation}
Furthermore, $\mathfrak{U}(\mathfrak{g})$ is
super-cocommutative, i.e., $T^s\circ\bigtriangleup=\bigtriangleup$,
where
\begin{equation}\label{T^s}
T^s(X\otimes Y)= (-1)^{p(X)p(Y)} Y\otimes X.
\end{equation}
Using (\ref{tensorprodukt}) and (\ref{T^s}) we get $f_1\cdot
f_2=(-1)^{p(f_1)p(f_2)}f_2\cdot f_1$. Hence, the sheaf
$\widehat{\mathcal{O}}_G$ is a sheaf of commutative associative
superalgebras with unit.

Further, $\bigwedge(\mathfrak{g}_{\bar 1})$ is also a cosuperalgebra
with comultiplication defined by
$$
\triangle_{\mathfrak{g}_{\bar 1}}(X)=X\otimes 1+ 1\otimes X,\,\,\,
\triangle_{\mathfrak{g}_{\bar 1}}(X_1\wedge \cdots \wedge X_r)=
\triangle_{\mathfrak{g}_{\bar 1}}(X_1)\wedge \cdots \wedge
\triangle_{\mathfrak{g}_{\bar 1}}(X_r),
$$
where $X,X_i\in \mathfrak{g}_{\bar 1}$. As above, this permits to
regard $\Hom(\bigwedge(\mathfrak{g}_{\bar 1}),\mathcal{F}_G)$ as a
sheaf of superalgebras which we may identify with the sheaf of
superalgebras $\mathcal{F}_G\otimes \bigwedge(\mathfrak{g}^*_{\bar
1})$. Moreover, the homomorphism $\gamma$ given by
(\ref{isomorphism}) is a homomorphism of
 cosuperalgebras. It follows that the mapping
$\widehat{\mathcal{O}}_G\to
\Hom (\bigwedge(\mathfrak{g}_{\bar 1}),\mathcal{F}_G)$, given by
$f\mapsto f\circ \gamma$, is an isomorphism of
sheaves of superalgebras.
%
%
Hence, $\widehat{\mathcal{O}}_G\simeq \mathcal{F}_G\otimes
\bigwedge(\mathfrak{g}^*_{\bar 1})$, and
$(G,\widehat{\mathcal{O}}_G)$ is a supermanifold. Clearly, it is
split and corresponds to the trivial bundle over $G$ with the fibre
$\bigwedge(\mathfrak{g}^*_{\bar 1})$.

Now we are able to define a structure of a Lie supergroup on
$(G,\widehat{\mathcal{O}}_G)$. The following formulas define the
multiplication morphism, the inversion morphism and the identity
morphism respectively (see \cite{BagSta}):
\begin{equation}\label{umnozh}
\begin{split}
\mu^*(f)(X\otimes Y)(g,h)&=f(X\cdot \al_G(g)(Y))(gh);\\
\iota^*(f)(X)(g)&=f(\al_G(g^{-1})(S(X)))(g^{-1});\\
\varepsilon^*(f)&=f(1)(e).
\end{split}
\end{equation}
Here $X,Y\in\mathfrak{U}(\mathfrak{g}),\, f\in
\widehat{\mathcal{O}}_G,\, g,\,h\in G$, and we identify the
enveloping superalgebra
$\mathfrak{U}(\mathfrak{g}\oplus\mathfrak{g})$ with the tensor
product $\mathfrak{U}(\mathfrak{g})\otimes
\mathfrak{U}(\mathfrak{g})$. The group axioms can be easily
verified, using the Hopf (super)algebra axioms. Note that
$(G,\widehat{\mathcal{O}}_G)$ corresponds to the Harish-Chandra pair
$(G,\mathfrak{g})$ and $\alpha_G=\Ad_G$.


Let $(\Phi, \varphi)$ be a morphism of Harish-Chandra pairs
$(G,\mathfrak{g})\to (H,\mathfrak{h})$. Then we can define a
morphism $F((\Phi, \varphi)) = \Psi: (G,\widehat{\mathcal{O}}_G) \to
(H,\widehat{\mathcal{O}}_H)$ by the following formula:
\begin{equation}\label{mor}
\Psi_{\red}= \Phi, \,\,\, \Psi^*(f)(X)(g)=f(\varphi(X))(\Phi(g)),
\,\,f\in \widehat{\mathcal{O}}_H, \,X\in
\mathfrak{U}(\mathfrak{g}),\, g\in G.
\end{equation}
Let us prove that $\Psi$ is a homomorphism of Lie supergroups. We
should show that $\Psi\circ \mu_G=\mu_H \circ (\Psi\times \Psi)$,
where $\mu_G$ and $\mu_H$ are the multiplication morphisms of
$(G,\widehat{\mathcal{O}}_G)$ and $(H,\widehat{\mathcal{O}}_H)$
respectively. By definition, we have
$$
\begin{array}{cc}
\mu_G^*\circ \Psi^*(f)(X\otimes Y)(g,h)=
\Psi^*(f)(\alpha_G(h^{-1})(X)\cdot Y)(gh)=
\\
f(\varphi(\alpha_G(h^{-1})(X)\cdot Y)) (\Phi(gh))
\end{array}
$$
and
$$
\begin{array}{c}
 (\Psi^*\times \Psi^*)\circ \mu_H^* (f)(X\otimes Y)(g,h)=
\mu_H^* (f)(\varphi(X)\otimes\varphi(Y))(\Phi(g),\Phi(h))=\\
f(\alpha_H(\Phi(h^{-1}))(\varphi(X))\cdot
\varphi(Y))(\Phi(g)\Phi(h)).
\end{array}
$$
Now our assertion follows from the definition of a morphism of
Harish-Chandra pairs.

\medskip

\noindent {\it 2.2  Isomorphisms of objects.} Let
$(G,\mathcal{O}_G)$ be a Lie supergroup and $\mathfrak{g}$ the
corresponding Lie superalgebra. We want to prove that
$(G,\mathcal{O}_G)$ is isomorphic to the Lie supergroup $F(P)$ which
corresponds to the Harish-Chandra pair $P = (G,\mathfrak{g})$.
Actually, we are going to prove a more general assertion, and
therefore we first extend the functor $F$ to a wider class of
objects.

Let $H$ be a closed Lie subgroup of $G$. As above, putting
$$
\widehat{\mathcal{O}}_{G/H}(U) =
\Hom_{\mathfrak{U}(\mathfrak{g}_{\bar 0})}(\mathfrak{U}(\mathfrak{g}),
\mathcal{F}_{G/H}(U))
$$
for every open $U\subset G/H$, we get a sheaf of superalgebras
$\widehat{\mathcal{O}}_{G/H}$. By the same argument as above,
$(G/H,\widehat{\mathcal{O}}_{G/H})$ is a split supermanifold and
$\widehat{\mathcal{O}}_{G/H}$ is isomorphic to
$\mathcal{F}_{G/H}\otimes \bigwedge(\mathfrak{g}^*_{\bar 1})$ (see
\cite{kosz}). The isomorphism $\widehat{\mathcal{O}}_{G/H} \to
\Hom(\bigwedge(\mathfrak{g}_{\bar 1}),\mathcal{F}_G)\simeq
\mathcal{F}_{G/H}\otimes \bigwedge(\mathfrak{g}^*_{\bar 1})$ is
again given by the formula $f\mapsto f\circ \gamma$, where $\gamma$
is defined by (\ref{isomorphism}).


Further, let $\nu: (G,\mathcal{O}_G)\times (M,\mathcal{O}_M) \to
(M,\mathcal{O}_M)$ be a transitive action. For simplicity we will
denote the vector field $(X\otimes \id)\circ \nu^*$ also by $X$.
Denote by $H$ the stabilizer of a certain point $x\in M$ by the
action $\nu_{\red}$. Our next aim is to define a morphism of
supermanifolds $(G/H,\widehat{\mathcal{O}}_{G/H})\to
(M,\mathcal{O}_M)$. We will use the natural correspondence
$X\mapsto\widehat{X}$ between even vector fields on
$(M,\mathcal{O}_M)$ and vector fields on $M$ which is completely
determined by the relation $\widehat{X}(f_{\red})= (X(f))_{\red}$
for all $f\in\mathcal{O}_{M}$. Let $f\in\mathcal{O}_{M}(U)$, where
$U$ is an open set in ${M}$. Denote by $\beta$ the natural
biholomorphic mapping  $G/H\to M$, $gH\mapsto gx$. Let us define the
linear mapping $\Phi_{G/H}(f):
\mathfrak{U}(\mathfrak{g})\to\mathcal{F}_{G/H}(\beta^{-1}(U))$ by

\begin{equation}
\label{def} \Phi_{G/H}(f)(X):= (-1)^{p(X)p(f)}\beta^*(X(f))_{\red},
\end{equation}
where $X\in\mathfrak{U}(\mathfrak{g})$ and $X$ and $f$ are
homogeneous. If $X\in \mathfrak{g}_{\bar 0}$, denote by
$\overline{X}$ the corresponding vector field on $G/H$. Note that
$\beta^*\circ \widehat{X}= \overline{X}\circ \beta^*$. The mapping
$\Phi_{G/H}(f)$ is a homomorphism of
$\mathfrak{U}(\mathfrak{g}_{\bar 0})$-modules. In fact, for any
$X_i\in \mathfrak{g}_{\bar 0}$, $Y_j\in \mathfrak{g}_{\bar 1}$ we
have
$$
\begin{array}{c}
 \Phi_{G/H}(f)(X_1\cdots X_r\cdot Y_1\cdots Y_q)=
  (-1)^{p(Y_1\cdots  Y_q) p(f)}\beta^*((X_1\cdots
 X_r\cdot Y_1\cdots  Y_q)(f))_{\red}\\=
(-1)^{p(Y_1\cdots  Y_q) p(f)}\beta^*(\widehat{X}_1\cdots
\widehat{X}_r)[ (Y_1\cdots Y_q(f))_{\red}]= \\
(-1)^{p(Y_1\cdots  Y_q) p(f)} (\overline{X}_1\cdots
\overline{X}_r)\beta^*[ (Y_1\cdots Y_q(f))_{\red}]= \\
(\overline{X}_1\cdots \overline{X}_r)(\Phi_{G/H}(f)(Y_1\cdots Y_q)).
\end{array}
$$

\medskip

\l\label{comutativn} {\it  $\Phi_{G/H}:
\mathcal{O}_{M}\to\widehat{\mathcal{O}}_{G/H}$ is a homomorphism of
sheaves of superalgebras and  $(\beta,
\Phi_{G/H}):(G,\widehat{\mathcal{O}}_{G/H})\to (M,\mathcal{O}_{M})$
is a morphism of supermanifolds.}

\medskip

 To prove Proposition \ref{comutativn},  we need the
following two lemmas.

\medskip

\lem\label{comutativn1} {\it Let $X_1,\ldots,X_r\in
\mathfrak{g}_{\bar 1}$, then
  \begin{equation} \label{delta}
\bigtriangleup(X_1\cdots X_r)= \sum\limits_{a+b=r}(-1)^{|\tau|}
X_{k_1}\cdots X_{k_a}\otimes X_{l_1}\cdots X_{l_b},
\end{equation}
where  $k_1<\cdots< k_a$, $l_1<\cdots< l_b$ and $\tau$ is the
permutation such that
$$
\tau (k_1,\ldots, k_a, l_1,\ldots, l_b)=(1,\ldots,r).
$$}

\medskip
\noindent{\it Proof.} For $r=1$ the formula is just the definition
of $\bigtriangleup$. Further, using induction, we get
$$
\begin{array}{c}
\bigtriangleup(X_1\cdots X_{r+1})= \bigtriangleup(X_1\cdots
X_r)\cdot\bigtriangleup(X_{r+1})=
\\
(\sum\limits_{a+b=r}(-1)^{|\tau|} X_{k_1}\cdots X_{k_a}\otimes
X_{l_1}\cdots X_{l_b})\cdot (X_{r+1}\otimes 1+ 1\otimes X_{r+1})\\
 =
\sum \limits_{a+b=r}(-1)^{|\tau|+b} X_{k_1}\cdots X_{k_a}\cdot
X_{r+1}\otimes
X_{l_1}\cdots X_{l_b} +\\
\sum\limits_{a+b=r}(-1)^{|\tau|} X_{k_1}\cdots X_{k_a}\otimes
X_{l_1}\cdots X_{l_b}\cdot X_{r+1} =\\
\sum\limits_{a'+b'=r+1}(-1)^{|\tau'|} X_{k_1}\cdots
X_{k_{a'}}\otimes X_{l_1}\cdots X_{l_{b'}},
\end{array}
$$
where $k_1<\cdots< k_{a'}$, $l_1<\cdots< l_{b'}$ and
$$
\tau' (k_1,\ldots, k_a', l_1,\ldots, l_b')=(1,\ldots,r+1). \Box
$$

\medskip

\lem\label{comutativn2} {\it Let $X_1,\ldots,X_r\in
\mathfrak{g}_{\bar 1}$ and $f_1,f_2\in \mathcal{O}_M$, then
\begin{equation}\label{diff}
(X_1\cdots X_r)(f_1f_2)=\!\!\! \sum_{a+b=r}(-1)^{|\tau|+p(f_1)b}
(X_{k_1}\cdots X_{k_a})(f_1) (X_{l_1}\cdots X_{l_b})(f_2),
\end{equation}
where  $k_1<\cdots< k_a$, $l_1<\cdots< l_b$ and $ \tau (k_1,\ldots,
k_a, l_1,\ldots, l_b)=(1,\ldots,r).
$
}

\noindent{\it Proof.} For $r=1$ the formula is simply the Leibniz
rule. Again, using induction, we get
$$
\begin{array}{c}
(X_1\cdots X_{r+1})(f_1f_2) =\\
X_1(\sum\limits_{a+b=r}(-1)^{|\tau|+p(f_1)b} (X_{k_1}\cdots
X_{k_a})(f_1) (X_{l_1}\cdots X_{l_b})(f_2))=\\
\rule{0pt}{5mm}\sum\limits_{a+b=r}(-1)^{|\tau|+p(f_1)b} (X_1\cdot
X_{k_1}\cdots X_{k_a})(f_1) (X_{l_1}\cdots X_{l_b})(f_2) +\\
\rule{0pt}{5mm}
 (-1)^{|\tau|+p(f_1)b+p(f_1)+a}
(X_{k_1}\cdots X_{k_a})(f_1)(X_1\cdot X_{l_1}\cdots X_{l_b})(f_2) =\\
\rule{0pt}{5mm} \sum\limits_{a'+b'=r+1}(-1)^{|\tau'|+p(f_1)b'}
(X_{k_1}\cdots X_{k_{a'}})(f_1) (X_{l_1}\cdots X_{l_{b'}})(f_2),
\end{array}
$$
where $k_1<\cdots< k_{a'}$, $l_1<\cdots< l_{b'}$ and
$$ \tau'
(k_1,\ldots, k_a', l_1,\ldots, l_b')=(1,\ldots,r+1). \Box $$

\noindent{\it Proof} of Proposition \ref{comutativn}. We should
check the equality
$$
(\Phi_{G/H}(f_1)\cdot \Phi_{G/H}(f_2))(X)= (\Phi_{G/H}(f_1f_2))(X)
$$
for $X\in \mathfrak{U}(\mathfrak{g})$, $f_1,f_2\in\mathcal{O}_{M}$.
Without loss of generality we may assume that $X=X_1\cdots X_r$,
$X_i\in \mathfrak{g}_{\bar 1}$, and that $f_1, f_2$ are homogeneous.
Using (\ref{delta}), we get
$$
\begin{array}{c}
(\Phi_{G/H}(f_1)\cdot \Phi_{G/H}(f_2))(X_1\cdots X_{r})=\\
\rule{0pt}{5mm} \operatorname{Mult}_{\mathcal{F}_G}
(\Phi_{G/H}(f_1)\otimes \Phi_{G/H}(f_2))
(\sum\limits_{a+b=r}(-1)^{|\tau|} X_{k_1}\cdots
X_{k_a}\otimes X_{l_1}\cdots X_{l_b})=\\
\rule{0pt}{5mm} \sum\limits_{a+b=r}(-1)^{|\tau|+p(f_2)a}
\Phi_{G/H}(f_1)(X_{k_1}\cdots
X_{k_a})\Phi_{G/H}(f_2)( X_{l_1}\cdots X_{l_b})=\\
\rule{0pt}{5mm}
\sum\limits_{a+b=r}(-1)^{|\tau|+p(f_2)a}(-1)^{p(f_1)a + p(f_2)b}
\beta^*{[}(X_{k_1}\cdots X_{k_a})(f_1) (X_{l_1}\cdots
X_{l_b})(f_2){]}_{\red}=\\
\rule{0pt}{5mm} \sum\limits_{a+b=r}(-1)^{|\tau|+p(f_2)r +
p(f_1)a}\beta^* {[} (X_{k_1}\cdots X_{k_a})(f_1) (X_{l_1}\cdots
X_{l_b})(f_2){]}_{\red}.
\end{array}
$$
On the other hand, by (\ref{diff}) we have
$$
\begin{array}{c}
(\Phi_{G/H}(f_1f_2))(X_1\cdots X_{r})=
(-1)^{r(p(f_1)+p(f_2))}\beta^*
 ((X_1\cdots X_{r})(f_1f_2))_{\red}=\\
\rule{0pt}{6mm}
(-1)^{r(p(f_1)+p(f_2))}\sum\limits_{a+b=r}(-1)^{|\tau|+p(f_1)b}
\beta^*{[}
(X_{k_1}\cdots X_{k_a})(f_1) (X_{l_1}\cdots X_{l_b})(f_2){]}_{\red}=\\
\rule{0pt}{5mm} \sum\limits_{a+b=r}(-1)^{|\tau|+p(f_2)r + p(f_1)a}
\beta^* {[}(X_{k_1}\cdots X_{k_a})(f_1) (X_{l_1}\cdots
X_{l_b})(f_2){]}_{\red}.
\end{array}
$$
The equality proves the first assertion of Proposition
\ref{comutativn}. The second assertion follows from the first one.
The proof is complete. $\Box$

\medskip

We may consider the special case when the odd dimension of
$(M,\mathcal{O}_M)$ is equal to the odd dimension of
$(G,\mathcal{O}_G)$. Later we will see that this condition is
equivalent to the following one:
$$
\dim (H,\mathcal{O}_{H})=\dim H| 0,
$$
where $(H,\mathcal{O}_{H})$ is the stabilizer of $x$ (see below). In
other words, $(H,\mathcal{O}_{H})=(H,\mathcal{F}_{H})$
 is an ordinary Lie subgroup of $G=(G,\mathcal{F}_{G})$.
\bigskip

\l\label{Phi -isomorph} {\it The morphism $(\beta, \Phi_{G/H}):
(G/H,\widehat{\mathcal{O}}_{G/H})\to (M,\mathcal{O}_{M})$ is a
submersion. If in addition $\dim (M,\mathcal{O}_M)=\dim M| \dim
\mathfrak{g}_{\bar 1}$, then $(\beta, \Phi_{G/H})$ is an isomorphism
of supermanifolds. In particular, all complex homogeneous
supermanifolds of this kind are split. }

\medskip
\noindent{\it Proof.} Let $y\in M$. Denote by $\mathfrak{m}_y$ and
$\widehat{\mathfrak{m}}_y$ the maximal ideals of the local
superalgebras $(\mathcal{O}_{M})_y$ and
$(\widehat{\mathcal{O}}_{G/H})_y$ respectively. It is easy to see
that
$$
\begin{array}{c}
\widehat{\mathfrak{m}}_y= \{ h\in
(\widehat{\mathcal{O}}_{G/H})_y\,\mid \,h(1)(y)=0 \},\\
\widehat{\mathfrak{m}}^2_y= \{ h\in \widehat{\mathfrak{m}}_y\,\mid
\,h(X)(y)=0\,\, \,\text{for all} \,X \in \mathfrak{g} \}.
\end{array}
$$
Note that
$\Phi_{G/H}(\mathfrak{m}_y)\subset\widehat{\mathfrak{m}}_y$. Let us
take $f\in \mathfrak{m}_y \setminus\mathfrak{m}_y^2$. The action
$\nu$ is transitive, hence there exists $X\in \mathfrak{g}$ such
that $(X(f))_{\red}(y)\ne 0$. Therefore, $\Phi_{G/H}(f)(X)(y)\ne 0$
and $\Phi_{G/H}(f) \in \widehat{\mathfrak{m}}_y \setminus
\widehat{\mathfrak{m}}_y^2$. It follows that the induced map
$\mathfrak{m}_y/ \mathfrak{m}_y^2\to \widehat{\mathfrak{m}}_y/
\widehat{\mathfrak{m}}_y^2$ is injective. Hence, the dual map
$(\widehat{\mathfrak{m}}_y/ \widehat{\mathfrak{m}}_y^2)^*\to
(\mathfrak{m}_y/ \mathfrak{m}_y^2)^*$ (or the differential of
$(\beta,\Phi_{G/H})$ at the point $y$) is surjective. Hence $(\beta,
\Phi_{G/H})$ is a submersion. Further, since
$$
\dim (M,\mathcal{O}_{M})=\dim (G/H,\widehat{\mathcal{O}}_{G/H})=\dim
M | \dim \mathfrak{g}_{\bar 1},
$$
we get that the differential is an isomorphism at every point $y\in
M$. Hence, $(\beta,\Phi_{G/H})$ is a local isomorphism (see
\cite{ley}, Inverse Function Theorem). But the mapping $\beta$ is
bijective, hence $(\beta,\Phi_{G/H})$ is an isomorphism. $\Box$

\medskip

In the case when $(M,\mathcal{O}_{M}) = (G,\mathcal{O}_{G})$ and
$\nu = \mu$ we get

\medskip

\noindent {\bf Corollary}.\label{sledstvie} {\it All complex
supergroups are split supermanifolds.}

\medskip

This fact also follows from the results of \cite{Mol}. Note that not all
complex homogeneous supermanifolds are split. Some examples can be found in
\cite{oniosp}.

Now we return to the correspondence between Lie supergroups and
Harish-Chandra pairs. In the case when $(M,\mathcal{O}_{M}) =
(G,\mathcal{O}_{G})$ formula (\ref{def}) defines a homomorphism of
sheaves of superalgebras $\Phi_{G}:
\mathcal{O}_{G}\to\widehat{\mathcal{O}}_{G}$ if we put $x=e$ and
$\beta=\id$. Define by $\Phi_{G}\times \Phi_{G}$ the second
component of the morphism
$$
(\id,\Phi_{G})\times (\id,\Phi_{G}):  (G,
\widehat{\mathcal{O}}_G)\times (G, \widehat{\mathcal{O}}_G) \to (G,
\mathcal{O}_G)\times (G, \mathcal{O}_G).
$$

\medskip

\lem\label{ravenstvo} {\it We have $\Phi_{G}\times
\Phi_{G}=\Phi_{G\times G}$.}

\medskip
\noindent{\it Proof.} It is sufficient to check the equality
$$
(\Phi_{G}\times \Phi_{G})\mid_{\pr^*_i(\mathcal{O}_G)}=\Phi_{G\times
G}\mid_{\pr^*_i(\mathcal{O}_G)},\,\, i=1,2,
$$
 where $\pr_i:(G,\mathcal{O}_G)\times (G,\mathcal{O}_G)\to
(G,\mathcal{O}_G)$ is the projection onto the $i$-th factor. Let
$\widehat{\pr}_i:(G,\widehat{\mathcal{O}}_G)\times
(G,\widehat{\mathcal{O}}_G)\to (G,\widehat{\mathcal{O}}_G)$ be also
the projection onto the $i$-th factor and $h\in
\widehat{\mathcal{O}}_G$. For example $\widehat{\pr}^*_1(h)$ has the
following form as a $\mathfrak{U}((\mathfrak{g}\oplus
\mathfrak{g})_{\bar 0})$-module homomorphism of
$\mathfrak{U}(\mathfrak{g}\oplus \mathfrak{g}) \simeq
\mathfrak{U}(\mathfrak{g})\otimes \mathfrak{U}(\mathfrak{g})$ to
$\mathcal{F}_{G\times G}$:
\begin{equation}\label{case}
\widehat{\pr}^*_1(h)(X^r\cdot Y^q)(g_1,g_2)= \left\{
                                                                 \begin{array}{ll}
                                                                   0, & \hbox{if}\,q\ne 0; \\
                                                                   h(X^q)(g_1), &
\hbox{if}\,q=0.
                                                                 \end{array}
                                                               \right.
\end{equation}
Here $X^r:= X_1\cdots X_r$, $Y^q:= Y_1\cdots Y_q$, where $X_i$ are
from the first copy of $\mathfrak{g}$ and $Y_j$ are from the second
one.

Let us take $f\in (\mathcal{O}_G)_{\bar i}$. By definition of
$\Phi_G\times \Phi_G$ and by (\ref{case}) we get:
$$
\begin{array}{rl}
(\Phi_G\times \Phi_G)(\pr_1^*(f))(X^r\cdot Y^q)(g_1,g_2)&= \\
\widehat{\pr}^*_1(\Phi_G(f))(X^r\cdot Y^q)(g_1,g_2)&= \left\{
\begin{array}{ll}
                                                                   0, & \hbox{if}\,q\ne 0; \\
                                                                  \Phi_G(f)(X^r)(g_1), &
\hbox{if}\,q=0.
                                                                 \end{array}
                                                               \right.
\end{array}
$$
On the other hand,
$$
\begin{array}{rl}
\Phi_{G\times G}(\pr^*_1(f)) (X^r\cdot Y^q)(g_1,g_2)=&
(-1)^{p(X^r\cdot Y^q)p(f)}[X^r\cdot Y^q
(\pr^*_1(f))]_{\red}(g_1,g_2)\\\rule{0pt}{7mm} =&\left\{
 \begin{array}{ll}
  0, & \hbox{if}\,q\ne 0; \\
  (-1)^{p(X^r) p(f)}(X^r (f))_{\red}(g_1), & \hbox{if}\,q=0.
   \end{array}
     \right.
\end{array}
$$
This completes the proof. $\Box$

\bigskip

\l\label{Phi -isomorph} {\it $(\id,\Phi_G): (G,\widehat{\mathcal{O}}_G)\to
(G,\mathcal{O}_G)$ is an isomorphism of Lie supergroups.}

\medskip
\noindent{\it Proof.} Due to Lemma \ref{ravenstvo}, we should check
that $(\Phi_{G\times G}\circ\mu^*)(f) = ((\hat\mu)^*\circ\Phi_G)(f)$
for all $f\in \mathcal{O}_G$, where $\hat\mu$ is the multiplication
morphism for $(G,\widehat{\mathcal{O}}_G)$. Let $X^r$ and $Y^q$ be
as in the proof of Lemma \ref{ravenstvo}. Recall that in
Preliminaries the morphism $\delta_x:(\pt,\mathbb{K})\to
(M,\mathcal{O}_M)$ was defined. Obviously we have
\begin{equation}\label{Phi_G times G}
\begin{split}
\Phi_{G\times G}(\mu^*(f))(X^r\cdot Y^q)(g_1,g_2) = (-1)^{p(X^r\cdot
Y^q)p(f)} (\delta_{g_1}^*\otimes \delta_{g_2}^*) \circ X^r\cdot Y^q
\circ\\
\mu^*(f)= (-1)^{p(X^r\cdot Y^q)p(f)}(\delta_{g_1}^*\circ X^r\otimes
\delta_{g_2}^*\circ Y^q)\circ \mu^*(f), \,\,g_1,g_2 \in G.
\end{split}
\end{equation}
We will use the following equalities:
\begin{equation}\label{ravenstva raznye}
\begin{split}
(r_{g_1}^*\otimes \id )\circ \mu^*&= (\id \otimes\, r_{g_1}^*\circ
\omega_{g_1}^*) \circ \mu^*,\\
r_{g_i}^* \circ X&=X \circ r_{g_i}^*, \,\,i=1,2;\\
\delta_{g_2}^*\circ r_{g_1}^* = & \,\,\delta_{g_2g_1}^*,\,\,\,
\delta_{g_2g_1}^*\circ \omega_{g_1}^*=\delta_{g_1g_2}^*.
\end{split}
\end{equation}
Here  $g_1,\,g_2\in G$, $X\in \mathfrak{U}(\mathfrak{g})$. By
(\ref{ravenstva raznye}) we get
\begin{equation}\label{prodolzhenie}
\begin{split}
(\delta_{g_1}^*\circ X^r\otimes & \delta_{g_2}^*\circ Y^q) \circ
\mu^*(f)= (\delta_{e}^*\circ r_{g_1}^*\circ X^r\,\otimes\,\,
\delta_{g_2}^*\circ
Y^q)\circ \mu^*(f)=\\
&(\delta_{e}^*\circ X^r\otimes\, \delta_{g_2}^*\circ Y^q)\circ
(r_{g_1}^* \otimes \,\id)\circ
 \mu^*(f)=\\
 & (\delta_{e}^*\circ X^r \otimes \delta_{g_2}^*\circ Y^q \circ r_{g_1}^*\circ \omega_{g_1}^*) \circ \mu^*(f)=\\
&(\delta_{e}^* \circ X^r\otimes \delta_{g_2}^*\circ r_{g_1}^*\circ
Y^q\circ \omega_{g_1}^*) \circ
\mu^*(f)=\\
&(\delta_{e}^* \circ X^r\otimes
\delta_{g_2g_1}^*\circ Y^q \circ \omega_{g_1}^*)\circ \mu^*(f)=\\
&(\delta_{e}^* \circ X^r\otimes \delta_{g_1g_2}^*\circ
\omega_{g_1^{-1}}^*\circ Y^q\circ \omega_{g_1}^*)\circ
\mu^*(f)=\\
&(\delta_{e}^*\circ  X^r \otimes \delta_{g_1g_2}^* \circ\Ad_G(g_1)(
Y^q))\circ \mu^*(f).
\end{split}
\end{equation}

\noindent By induction it is easy to check that
\begin{equation}\label{mu}
\begin{split}
&(Y^q)(f)= (-1)^{A(Y^q)}(\delta_{e}^*\circ Y_q\otimes \cdots \otimes
\delta_{e}^*\circ Y_1\otimes \id)\circ(\mu^q)^*(f),\,\,\text{where}\\
&A(Y^q)=p(Y_{q-1})p(Y_q)+ p(Y_{q-2})p(Y_{q-1}\cdot Y_q)+\cdots+
p(Y_1)p(Y_{2}\cdots Y_q).
\end{split}
\end{equation}
Here $\mu^q$ is the multiplication morphism of $q+1$ copies of
$(G,\mathcal{O}_G)$. Indeed, for $q=1$ the assertion (\ref{mu}) is
just the definition of a right invariant vector field. Further,
$$
\begin{array}{c}
(Y^q)(f)= Y_1((-1)^{A(Y_{2}\cdots Y_q)}(\delta_{e}^*\circ Y_q\otimes
\cdots
\otimes \,\delta_{e}^*\circ Y_2\otimes \id)(\mu^{q-1})^*(f))=\\
(-1)^{A(Y_{2}\cdots Y_q)}   (-1)^{p(Y_1)p(Y_{2}\cdots Y_q)}
(\delta_{e}^*\circ Y_q\otimes\cdots \otimes\, \delta_{e}^*\circ
Y_1\otimes
\id)(\mu^{q})^*(f)=\\
(-1)^{A(Y^q)}(\delta_{e}^*\circ Y_q\otimes\cdots \otimes\,
\delta_{e}^*\circ Y_1\otimes \id)(\mu^{q})^*(f).
\end{array}
$$

\noindent By (\ref{mu})
 we have
$$
\begin{array}{c}
 (\delta_{e}^*\circ X^r\otimes \delta_{g_1g_2}^*\circ\Ad_G(g_1)( Y^q))\circ \mu^*(f)=
(-1)^{A(X^r)+ A(Y^q)} (\delta_{e}^*\circ X_r \otimes
\cdots\\
\otimes\delta_{e}^*\circ  X_1 \otimes \delta_{e}^*\circ \Ad_G(g_1)(
Y_q)\otimes \cdots \otimes \delta_{e}^*\circ \Ad_G(g_1)( Y_1)\otimes
\delta_{g_1g_2}^*)
(\mu^{r+q})^*(f)=\\
(-1)^{A(X^r)+ A(Y^q)+ A(Y^q\cdot X^r)}(\delta_{g_1g_2}^*\circ
\Ad_G(g_1)( Y^q) \circ X^r) (f)=\\
(-1)^{A(X^r)+ A(Y^q)+ A(Y^q\cdot
X^r)}(-1)^{p(X^r)p(Y^q)}(\delta_{g_1g_2}^*\circ X^r\circ \Ad_G(g_1)(
Y^q)) (f)=\\
(\delta_{g_1g_2}^*\circ X^r\circ \Ad_G(g_1)( Y^q)) (f).
\end{array}
$$
On the other hand,
$$
\begin{array}{c}
\mu^*(\Phi(f))(X^r\cdot Y^q)(g_1,g_2)=
\Phi(f)(X^r\cdot\Ad_G(g_1) (Y^q))(g_1g_2)=\\
(-1)^{p(X^r\cdot Y^q)p(f)}(\delta_{g_1g_2}^*\circ X^r \cdot
\Ad_G(g_1) (Y^q)) (f).
\end{array}
$$
This completes the proof. $\Box$

\medskip

\noindent{\it 2.3 The bijection between morphisms}. Let
$(G,\mathfrak{g})$ and $(H,\mathfrak{h})$ be two Harish-Chandra
pairs and $(G,\widehat{\mathcal{O}}_G)$ and
$(H,\widehat{\mathcal{O}}_H)$ be the corresponding Lie supergroups
with multiplication morphisms $\mu_G$ and $\mu_H$ respectively. Let
$\Psi: (G,\widehat{\mathcal{O}}_G)\to (H,\widehat{\mathcal{O}}_H)$
be a homomorphism of Lie supergroups. Let $X_e\in \mathfrak{g}$ and
$X=(X_e\otimes \id)\circ \mu^*_G$ be the corresponding right
invariant vector field on $(G,\widehat{\mathcal{O}}_G)$ and $Y=
((\d\Psi)_e X_e\otimes \id)\circ \mu^*_H$. Then the vector fields
$X$ and $Y$ are $\Psi$-related, i.e.
$$
X(\Psi^*(f))= \Psi^*(Y(f)),\,\,\,\, f\in \widehat{\mathcal{O}}_H.
$$
Now we are able to prove that $\Psi$ depends only on $\Psi_{\red}$
and $(\d \Psi)_e$. Indeed,
$$
\begin{array}{c} [\Psi^*(f)(X_e)](g)= (-1)^{p(X)p(f)}
[X(\Psi^*(f))]_{\red}(g)=\\ (-1)^{p(X)p(f)}[\Psi^*(Y(f))]_{\red}(g)=
(-1)^{p(X)p(f)} [Y(f)]_{\red}(\Psi_{\red}(g))=\\
 f((\d\Psi)_e X_e) (\Psi_{\red}(g)),
\end{array}
$$
 where $X\in \mathfrak{g}$, $f\in \widehat{\mathcal{O}}_H$, $g\in
 G$.
It follows that all homomorphisms of $(G,\widehat{\mathcal{O}}_G)$
to $(H,\widehat{\mathcal{O}}_H)$ have the form (\ref{mor}) if we put
$\varphi=(\d\Psi)_e$, $\Phi=\Psi_{\red}$. Hence the map
$F:\Hom(X,Y)\to \Hom(F(X), F(Y))$ is surjective. The injectivity of
the map $F:\Hom(X,Y)\to \Hom(F(X), F(Y))$ is obvious.
\medskip

\noindent{\it 2.4 The main result. }We have proved the following
theorem.

\medskip
\t\label{main} {\it The category of complex Lie supergroups is
equivalent to the category of complex Harish-Chandra pairs.  }

\medskip

Theorem \ref{main} implies some important consequences: the
existence of a Lie supergroup for a given Lie superalgebra, the
existence of a Lie subsupergroup for a given Lie subsuperalgebra.
(The last assertion we will discuss below.) Using Theorem
\ref{main}, these two assertions can be proven in the
complex-analytic case as in \cite{Kostant}, Corollary to Theorem
$3.7$ and Theorem $3.8$. 

\medskip

\noindent {\bf Remark.} The same method can be used to prove Theorem
\ref{main} in the category of affine algebraic supergroups in the
sense of \cite{Westra}.

\bigskip


\begin{center}
{\bf 3. Homogeneous supermanifolds}
\end{center}

\medskip

Suppose that a closed Lie subsupergroup $(H,\mathcal{O}_H)$ of
$(G,\mathcal{O}_G)$ (this means that the Lie subgroup $H$ is closed
in $G$) is given. Consider the corresponding coset superspace
$(G/H,\mathcal{O}_{G/H})$ (see \cite{Fio_Varad, Kostant}). Denote by
$\mu_{G\times H}$ the composition of the morphisms
$$
(G,\mathcal{O}_G)\times (H,\mathcal{O}_H)\hookrightarrow
(G,\mathcal{O}_G)\times (G,\mathcal{O}_G)\xrightarrow{\mu}
(G,\mathcal{O}_G),
$$
 by $\pr_1:(G,\mathcal{O}_G)\times (H,\mathcal{O}_H)\to
(G,\mathcal{O}_G)$ the projection onto the first factor and by $\pi$
the natural mapping $G\to G/H$, $g\mapsto gH$. Let us take $U\subset
G/H$ open.  Then
\begin{equation}\label{sheaf O_G/H}
\mathcal{O}_{G/H}(U)= \{ f\in \mathcal{O}_G(\pi^{-1}(U))\mid
(\mu_{G\times H})^*(f)=\pr^*_1(f)\}.
\end{equation}
Denote by $\nu: (G,\mathcal{O}_G)\times (G/H,\mathcal{O}_{G/H})\to
(G/H,\mathcal{O}_{G/H})$ the natural supergroup action. It is given
by $\nu^*(f)=\mu^*(f)$, where $f\in \mathcal{O}_{G/H}(U)\subset
\mathcal{O}_{G}(\pi^{-1}(U))$.

A Harish-Chandra pair $(H,\mathfrak{h})$ is called a {\it
Harish-Chandra subpair} of a Harish-Chandra pair $(G,\mathfrak{g})$
if $H$ is a Lie subgroup of $G$ and $\mathfrak{h}$ is a  Lie
subsuperalgebra of $\mathfrak{g}$, s.t. $\mathfrak{h}_{\bar 0}=\Lie
H$ and $\alpha_H=\alpha_G|H$. There is a correspondence between
Harish-Chandra subpairs of $(G,\mathfrak{g})$ and Lie subsupergroups
of $(G,\mathcal O_G)$. (The Lie supergroup $(G,\mathcal O_G)$
corresponds to the Harish-Chandra pair $(G,\mathfrak{g})$.) More
precisely, let $\varphi:(H,\mathcal O_H) \to (G,\mathcal O_G)$ be a
Lie subsupergroup. Then the corresponding Harish-Chandra pair
$(H,\mathfrak{h})$ is a Harish-Chandra subpair, because $H\subset G$
is a Lie subgroup and $\mathfrak{h}=(\d \varphi)_e \subset
\mathfrak{g}$
 is a Lie subalgebra, s.t. $\mathfrak{h}_{\bar 0}=\Lie
H$ and $\alpha_H=\alpha_G|H$. Further, let $\varphi:(H,\mathcal O_H)
\to (G,\mathcal O_G)$ and $\varphi':(H',\mathcal O_{H'}) \to
(G,\mathcal O_G)$ be two Lie subsupergroups which determine the same
Harish-Chandra pair $(H,\mathfrak{h})$. We claim that there is an
isomorphism of Lie supergroups $\psi:(H,\mathcal O_H) \to
(H',\mathcal O_{H'})$, such that $\varphi'\circ \psi=\varphi$. As we
have seen above, any homomorphism of Lie supergroups is determined
by its underlying map and its differential at the point $e$. To
define $\psi$, we put $\psi_{\red}=\id:H\to H'$ (note that  $H=H'$)
and $(\d\psi)_e = (\d\varphi')_e^{-1}\circ (\d\varphi)_e$.

Conversely, let $(H,\mathfrak{h})$ be a Harish-Chandra subpair of
$(G,\mathfrak{g})$. Then we get the Lie supergroup $(H,\mathcal
O_H)$ using the construction from $2.1$. There is a natural
homomorphism $\varphi: (H,\mathcal O_H) \to (G,\mathcal O_G)$, where
$\varphi_{\red}:H\to G$ is the inclusion and $\varphi^*:\mathcal O_G
\to \mathcal O_H$ is given by
$$
\varphi^*(f)(X)(h) = f(X)(\varphi_{\red}(h)),\,\,\, X\in
\mathfrak{U}(\mathfrak{h})\subset\mathfrak{U}(\mathfrak{g}),\,\,\,\,
h\in H.
$$
Clearly, the Harish-Chandra subpair which corresponds to the Lie
subsupergroup $(H,\mathcal{O}_H)$ coincides with $(H,\mathfrak{h})$.

Let $\nu: (G,\mathcal{O}_G)\times (M,\mathcal{O}_M)\to
(M,\mathcal{O}_M)$ be a transitive action. Denote by $\nu_{x}$,
where $x\in M$,  the composition of the morphisms
$$
(G,\mathcal{O}_G)=(G,\mathcal{O}_G)\times
(x,\mathbb{K})\stackrel{\id\times \delta_x}{\longrightarrow}
(G,\mathcal{O}_G)\times (M,\mathcal{O}_M)\xrightarrow{\nu}
(M,\mathcal{O}_M).
$$

\medskip

\lem\label{dlya glupyh} {\it We have $\ev_x\circ
\overline{\nu}(X)=(\d\nu_{x})_{e}(X_e)$, $X\in \mathfrak{g}$.}

\medskip

\noindent{\it Proof.} By definition we get
$$
\ev_x(\overline{\nu}(X))(f)=[\overline{\nu}(X)(f)]_{\red}(x),\,\,\,\,\,\,\,
(\d\nu_{x})_{e}(X_e)(f)=X_e\circ \nu^*_x(f)
$$
 for $f\in
(\mathcal O_M)_x$. Let $\delta_x(h):=(h)_{\red}(x)$, $h\in
\mathcal{O}_M$. Then,
$$
\begin{array}{c}
\ev_x(\overline{\nu}(X))(f)= [(X_e\otimes \id)\circ \nu^*(f)]_{\red}
(x) = (X_e\otimes \delta_x)\circ \nu^*(f) = \\
X_e\circ (\id \otimes \delta_x)\circ \nu^*(f)= X_e \circ \nu^*_x(f)
= (\d\nu_{x})_{e}(X_e)(f).\,\,\Box
\end{array}
$$

\noindent{\it Remark.} By the axioms of an action we have
$\nu_x=\nu_{gx}\circ r_{g^{-1}}$ for all $g\in G$. Using Lemma
\ref{dlya glupyh} we get that a supermanifold $(M,\mathcal O_M)$ is
$(G,\mathcal O_G)$-homogeneous if and only if $\nu_{\red}$ is a
transitive action of $G$ on $M$ and $(\d \nu_x)_e$ is surjective for
some $x\in M$.


 As in
\cite{Kostant}, we can define the stationary subsupergroup
$(G_x,\mathcal{O}_{G_x})$ of the point $x$ in the following way.
Consider the Harish-Chandra subpair $(G_x,\mathfrak{g}_x)$ of
$(G,\mathfrak{g})$, where $G_x$ is the stabilizer of $x$ and
$\mathfrak{g}_x=\Ker (\d \nu_x)_e$. A subsupergroup
$(G_x,\mathcal{O}_{G_x})$ is called the {\it stabilizer of} $x$, if
it determines $(G_x,\mathfrak{g}_x)$. Further, assume that the
action $\nu$ is transitive. In this case in \cite{V} another
definition of the stabilizer of $x$ was given. It is easy to prove
that these two definitions are equivalent. Moreover,
$(M,\mathcal{O}_M)\simeq (G/G_x,\mathcal{O}_{G/G_x})$ and hence
$\dim (M,\mathcal{O}_M)= \dim (G,\mathcal{O}_G)- \dim
(G_x,\mathcal{O}_{G_x})$ (see, \cite{Kostant,V}).

\bigskip

\begin{center}
{\bf 4. Homogeneous split supermanifolds}
\end{center}

\medskip

In this section we will consider only complex supermanifolds. Note
that all real super\-manifolds are split.

Let us introduce a new category $\verb"SSM"$ (split supermanifolds).
We put
$$
\begin{array}{rl}
\operatorname{Ob}\, \verb"SSM"=\{ (M,\bigwedge
\mathcal{E}_M)\,\,|\,\, \mathcal{E}_M\,\,\text{ is a locally free
sheaf on}\,\, M\}.
\end{array}
$$
Equivalently, we can say that $\operatorname{Ob}\, \verb"SSM"$
consists of all split supermanifolds $(M,\mathcal{O}_M)$ with a
fixed isomorphism $\mathcal{O}_M\simeq \bigwedge \mathcal{E}_M$ for
a certain locally free sheaf $\mathcal{E}_M$ on $M$. Note that
$\mathcal{O}_M$ is naturally $\mathbb{Z}$-graded by
$(\mathcal{O}_M)_p\simeq \bigwedge^p \mathcal{E}_M$. All the time we
will consider this $\mathbb{Z}$-grading. Further, if $X,Y\in
\operatorname{Ob}\, \verb"SSM"$ we put
$$
\begin{array}{rl}
\Hom(X,Y)=&\text{all morphisms of $X$ to $Y$,}\\
 &\text{preserving the
$\mathbb{Z}$-gradings}.
\end{array}
$$

As in the category of supermanifolds, we can define in $\verb"SSM"$
a group object (split Lie supergroup), an action  of a split Lie
supergroup on a split supermanifold (split action) and a homogeneous
split supermanifold. More precisely, we get these notions if we
consider in the definitions \ref{def Liesgroup}, \ref{def sgroup
action}, \ref{def homogen supermnf} morphisms and objects from
$\verb"SSM"$.

Let $(M,\bigwedge \mathcal{E}_M)$ and $(N,\bigwedge \mathcal{E}_N)$
be two split supermanifolds, where $\mathcal{E}_M$ and
$\mathcal{E}_N$ are the sheaves of sections of vector bundles $\bold
E_M$ and $\bold E_N$ respectively. The direct product in the
category $\verb"SSM"$ is defined by:
$$
(M,\bigwedge \mathcal{E}_M)\times (N,\bigwedge \mathcal{E}_N):=
(M\times N,\bigwedge(\mathcal{E}_M\oplus\mathcal{E}_N)).
$$
Here the fixed $\mathbb{Z}$-grading is given by
$$
\bigwedge^p(\mathcal{E}_M\oplus\mathcal{E}_N)=\bigoplus_{t+s=p}
\bigwedge^r \mathcal{E}_M\otimes \bigwedge^s \mathcal{E}_N.
$$
It is easy to see that this definition agrees with the definition of
the direct product in the category of supermanifolds, see
\cite{ley}.

There is a functor, say $\gr$, from the category of supermanifolds
to the category of split supermanifolds. Let us briefly describe
this construction (see, e.g., \cite{onipi,oniosp}). Let
$(M,\mathcal{O}_M)$ be a supermanifold. As above denote by
$\mathcal{J}_M\subset \mathcal{O}_M$ the subsheaf of ideals
generated by odd elements of $\mathcal{O}_M$. Then by definition
$\gr(M,\mathcal{O}_M)$ is equal to the split supermanifold $(M,
\gr\mathcal{O}_M)$, where
$$
\gr\mathcal{O}_M= \bigoplus_{p\geq 0} (\gr\mathcal{O}_M)_p,\quad
\mathcal{J}_M^0:=\mathcal{O}_M, \quad (\gr\mathcal{O}_M)_p=
\mathcal{J}_M^p/\mathcal{J}_M^{p+1}.
$$
In this case $(\gr\mathcal{O}_M)_1$ is a locally free sheaf and
there is a natural isomorphism of $\gr\mathcal{O}_M$ onto $\bigwedge
(\gr\mathcal{O}_M)_1$. If
$\psi=(\psi_{\red},\psi^*):(M,\mathcal{O}_M)\to (N,\mathcal{O}_N)$
is a morphism, then $\gr(\psi)=(\psi_{\red},\gr(\psi^*))$, where
$\gr(\psi^*):\gr \mathcal{O}_N \to \gr \mathcal{O}_M$ is defined by
$$
\gr(\psi^*)(f+\mathcal{J}_N^p): = \psi^*(f)+\mathcal{J}_M^p
\,\,\text{for}\,\, f\in (\mathcal{J}_N)^{p-1}.
$$
Recall that by definition every morphism $\psi$ of supermanifolds is
even and as consequence sends $\mathcal{J}_N^p$ into
$\mathcal{J}_M^p$.

Let $(G,\mathcal{O}_G)$ be a Lie supergroup with the group morphisms
$\mu$, $\iota$ and $\varepsilon$. Then it is easy to see that $\gr
(G,\mathcal{O}_G)$ is a split Lie supergroup with the group
morphisms $\gr(\mu)$, $\gr(\iota)$ and $\gr(\varepsilon)$.
Similarly, an action $\nu:(G,\mathcal{O}_G)\times
(M,\mathcal{O}_M)\to (M,\mathcal{O}_M)$ gives rise to the action
$\gr(\nu):\gr(G,\mathcal{O}_G)\times \gr(M,\mathcal{O}_M)\to
\gr(M,\mathcal{O}_M)$.

Obviously a split Lie supergroup is a Lie supergroup. Furthermore,
the following result holds:

\medskip
\l\label{split Lie supergroup} {\it A Lie supergroup
$(G,\mathcal{O}_G)$ with the Lie superalgebra $\mathfrak{g}$
possesses a structure of a split Lie supergroup if and only if
$[\mathfrak{g}_{\bar 1}, \mathfrak{g}_{\bar 1}] = 0$. Any Lie
subsupergroup of a split Lie supergroup possesses a structure of a
split Lie supergroup.

}
\medskip

\noindent{\it Proof.} Let $(G,\mathcal{O}_G)$ be a split Lie
supergroup, $\mathcal{O}_G=\bigoplus_p (\mathcal{O}_G)_p$ the fixed
$\mathbb{Z}$-grading and $\mu$ the multiplication morphism. 
Let us prove that $[\mathfrak{g}_{\bar 1},\mathfrak{g}_{\bar
1}]=\{0\}$. It is enough to check that $[X,Y](f)=0$ for $f\in
(\mathcal{O}_G)_0$ and $f\in (\mathcal{O}_G)_1$, where $X,Y\in
\mathfrak{g}_{\bar 1}$. By (\ref{left inv vect field}) we get
$$
\begin{array}{c}
[X,Y](f)=(X\otimes \id)\circ\mu^*\circ (Y\otimes
\id )\circ\mu^*(f)+\\
(Y\otimes \id)\circ\mu^*\circ (X\otimes \id)\circ\mu^*(f)=\\ -((
Y\otimes X\otimes \id)+(X\otimes Y \otimes \id))\circ (\mu^2)^*(f),
\end{array}
$$
where $\mu^2$ is the multiplication morphism of three copies of
$(G,\mathcal{O}_G)$. Note that by definition of a split Lie
supergroup $(\mu^2)^*(f)\in (\mathcal{O}_{G\times G\times G})_0$, if
$f\in (\mathcal{O}_G)_0$ and $(\mu^2)^*(f)\in (\mathcal{O}_{G\times
G\times G})_1$, if $f\in (\mathcal{O}_G)_1$. It follows that
$$
(X\otimes Y\otimes \id)((\mu^2)^*(f))=(Y\otimes X \otimes
\id)((\mu^2)^*(f))=0.
$$

Conversely, let $(G,\mathcal{O}_G)$ be a Lie supergroup and
$[\mathfrak{g}_{\bar 1}, \mathfrak{g}_{\bar 1}] = 0$. As we have
seen above, the sheaf
$\mathcal{O}_G=\Hom_{\mathfrak{U}(\mathfrak{g}_{\bar
0})}(\mathfrak{U}(\mathfrak{g}),\mathcal{F}_G)$ is a
$\mathbb{Z}$-graded sheaf. Recall that the $\mathbb{Z}$-grading is
induced by the mapping $\mathcal{O}_G\to
\Hom(\bigwedge(\mathfrak{g}_{\bar 1}),\mathcal{F}_G)$, $f\mapsto
f\circ \gamma$, where $\gamma$ is defined by (\ref{isomorphism}).
More precisely,
$$
f\in (\mathcal{O}_G)_p \Leftrightarrow f\circ \gamma (X_1\wedge
\cdots \wedge X_r)=0 \,\, \text{for}\,\, r\ne p,\,\, X_i\in
\mathfrak{g}_{\bar 1}.
$$
 It
follows that any Lie supergroup is contained in $\operatorname{Ob}\,
\verb"SSM"$.
 We want to prove that the structure
morphisms (\ref{umnozh}) of $(G,\mathcal{O}_G)$ preserve this
$\mathbb{Z}$-grading. Let us check that
$\mu^*((\mathcal{O}_G)_p)\subset (\mathcal{O}_{G\times G})_p$,
$p\geq 0$. Further, if $[\mathfrak{g}_{\bar 1}, \mathfrak{g}_{\bar
1}] = 0$, then $\gamma$ is a homomorphism of algebras (not only of
coalgebras), see (\ref{isomorphism}). It follows that
\begin{equation}\label{O_G grading}
f\in (\mathcal{O}_G)_p \Leftrightarrow f (X_1 \cdots
 X_r)=0 \,\, \text{for}\,\, r\ne p,\,\, X_i\in
\mathfrak{g}_{\bar 1}.
\end{equation}

Let $X_1,\ldots,X_r$ be from the first copy of $\mathfrak{g}_{\bar
1}$ and $Y_1,\ldots, Y_q$ be from the second one, $r+q\ne p$,
$g,h\in G$ and $f\in (\mathcal{O}_G)_p$, then
$$
\mu^*(f)(X_1\cdots X_r\cdot Y_1\cdots Y_q)(g,h)= f(X_1\cdots X_r
\cdot \alpha_G(g)(Y_1\cdots Y_q))(gh)=0.
$$
This implies that $\mu^*(f)\in (\mathcal{O}_{G\times G})_p$. For the
inversion morphism the proof is similar.

The second assertion is obvious. Indeed, if
$(H,\mathcal{O}_H)\subset (G,\mathcal{O}_G)$ is a Lie subsupergroup,
$\mathfrak{h} = \Lie (H,\mathcal{O}_H)$, then
$\mathfrak{h}\subset\mathfrak{g}$ is a subsuperalgebra. Hence,
$[\mathfrak{h}_{\bar 1}, \mathfrak{h}_{\bar 1}]=0$.$\Box$

\medskip

 \noindent {\bf Corollary.} {\it A split Lie supergroup is a semi-direct
product of a usual Lie group $G$ and the (unique) connected
supergroup of purely odd dimension.}

\medskip

\lem\label{split action} {\it Let $\nu:(G,\mathcal{O}_G)\times
(M,\mathcal{O}_M)\to (M,\mathcal{O}_M)$ be a transitive action of a
Lie supergroup $(G,\mathcal{O}_G)$ on a supermanifold
$(M,\mathcal{O}_M)$, then the action
$$
\gr\nu:(G,\gr\mathcal{O}_G)\times (M,\gr\mathcal{O}_M)\to
(M,\gr\mathcal{O}_M)
$$
is also transitive. In particular, if $(M,\mathcal{O}_M)$ is split
and homogeneous, then it always admits a transitive split action.

}
\medskip

\noindent{\it Proof.}  Since $\gr\nu_{\red}=\nu_{\red}$ it is enough
to show that $\d(\gr{\nu}_x)_e$ is surjective for some $x\in M$ (see
Remark after
 Lemma  \ref{dlya glupyh}). Since $\d(\gr\nu_x)_e = \d(\nu_x)_e$,
the proof is complete. The second assertion follows from the
isomorphism $(M,\mathcal{O}_M)\simeq (M,\gr\mathcal{O}_M)$.$\Box$

\medskip

Let $H$ be a closed Lie subgroup of a Lie group $G$, $E$ a complex
vector space and $\theta:H\to \GL(E)$ be a holomorphic
representation. Denote by $\mathcal{E}^{\theta}$ the sheaf of
sections of the homogeneous vector bundle $\bold E^{\theta}$ which
corresponds to $\theta$, i.e., the quotient space of the direct
product $G\times E$ by the following action of $H$:
$$
(g,v)\stackrel{h}{\mapsto} (gh^{-1}, \theta(h)v), \,\,\, g\in G,\,
h\in H,\, v\in E.
$$

Furthermore, let $\pi:G\to G/H$ be the natural projection and
$U\subset G/H$ open. There is an injective homomorphism of sheaves
$\Phi_{\theta}: \mathcal{E}^{\theta} \to \pi_*(\mathcal{F}_G\otimes
E)$ given by
$$
\begin{array}{rl}
\mathcal{E}^{\theta}(U)\ni s&\mapsto f_s\in
\mathcal{F}_G(\pi^{-1}(U))\otimes E,\\
f_s(g)&:=g^{-1}s(gH),\,\,g\in G.
\end{array}
$$
It is well known that
$$
\Phi_{\theta}(\mathcal{E}^{\theta}(U))= \{ f\in
\mathcal{F}_G(\pi^{-1}(U))\otimes E\,\mid\, \theta(h)f(gh)= f(g),\,
g\in G,\,h\in H\}.
$$
Note that $\bigwedge \mathcal{E}^{\theta}=
\mathcal{E}^{\wedge\theta}$ is also a homogeneous bundle. An easy
computation shows that $\Phi_{\wedge\theta}:
\bigwedge\mathcal{E}^{\theta} \to \pi_*(\mathcal{F}_G\otimes
\bigwedge E)$ is a homomorphism of the sheaves of superalgebras.

Let $V$ be  a vector space. Our aim is now to describe the
isomorphism of sheaves of superalgebras
$$
\Psi_V: \mathcal{F}_G \otimes \bigwedge V^*\to \Hom (\bigwedge
V,\mathcal{F}_G)
$$
mentioned in $2.1$.  Let $(\xi_i)$ be a basis of $V$, $(\xi_i^*)$
the dual basis of $V^*$, $h\in \mathcal{F}_G$ and $i_1<\cdots< i_k$,
$j_1<\cdots< j_r$. Then
$$
\begin{array}{rl}
\Psi_V(h \xi_{i_1}^*\wedge \cdots \wedge
\xi_{i_k}^*)&=(-1)^{k(k-1)/2}h f^{\xi_{i_1}\wedge \cdots \wedge
\xi_{i_k}},\,\,\text{where}\\
f^{\xi_{i_1}\wedge \cdots \wedge \xi_{i_k}}(\xi_{j_1}\wedge \cdots
\wedge \xi_{j_r})&=\left\{
                    \begin{array}{ll}
                      1, & \hbox{$(j_1,\cdots, j_r)=(i_1,\cdots, i_k)$;} \\
                      0, & \hbox{$(j_1,\cdots, j_r)\ne (i_1,\cdots, i_k)$.}
                    \end{array}
                  \right.
\end{array}
$$
The direct computation shows that it is a homomorphism of the
sheaves of superalgebras.

\medskip
\l\label{quptient of split Lie supergroups} {\it Let
$(G,\mathcal{O}_G)$ be a split Lie supergroup,
$(H,\mathcal{O}_H)\subset (G,\mathcal{O}_G)$ a closed Lie
subsupergroup and $(M,\mathcal{O}_M):=(G/H,\mathcal{O}_{G/H})$. The
$\mathbb{Z}$-grading $\mathcal{O}_G=\bigoplus_{p\geq
0}(\mathcal{O}_G)_p$, where $(\mathcal{O}_G)_p$ is determined by
(\ref{O_G grading}), induces the $\mathbb{Z}$-grading
$(\mathcal{O}_M)_p$ on the subsheaf $\mathcal{O}_M\subset
\mathcal{O}_G$. Moreover, $(\mathcal{O}_M)_1$ is a locally free
sheaf and $(\mathcal{O}_M)_p=\bigwedge^p(\mathcal{O}_M)_1$.

In particular,
 the coset supermanifold
$(M,\mathcal{O}_M)$ is split and the natural action of the Lie
supergroup $(G,\mathcal{O}_G)$ on $(M,\mathcal{O}_M)$ is split.

}
\medskip

\noindent{\it Proof.} The sheaf $\mathcal{O}_M$ was defined by
(\ref{sheaf O_G/H}). The supergroup $(G,\mathcal{O}_G)$ is a split
Lie supergroup, hence $(\mu_{G\times H})^*((\mathcal{O}_{G})_p)
\subset (\mathcal{O}_{G\times H})_p$. Furthermore, it is easy to see
that $\pr^*_1((\mathcal{O}_{G})_p)\subset (\mathcal{O}_{G\times
H})_p$. Hence the sheaf $\mathcal{O}_M$ is $\mathbb{Z}$-graded by
the subsheaves
$$
(\mathcal{O}_{M})_p= \{ f\in (\mathcal{O}_G)_p \mid (\mu_{G\times
H})^*(f)=\pr^*_1(f)\}.
$$
Consider the representation $\psi:H\to GL((\mathfrak{g}_{\bar
1}/\mathfrak{h}_{\bar 1})^*)$ defined by
$$
[\psi(h)(v)](X+\mathfrak{h}_{\bar 1}) = v
(\Ad_G(h^{-1})(X)+\mathfrak{h}_{\bar 1}) \,\,\text{ for}\,\,h\in
H,\,\, v\in (\mathfrak{g}_{\bar 1}/\mathfrak{h}_{\bar 1})^*,\,\,
X\in \mathfrak{g}_{\bar 1}.
$$
Our goal now is to show that $\mathcal{O}_{M}\simeq \bigwedge
\mathcal{E}^{\psi}$ as sheaves of $\mathbb{Z}$-graded algebras.
Denote by $\Gamma$ the isomorphism of sheaves
$\mathcal{O}_G=\widehat{\mathcal{O}}_G\to
\Hom(\bigwedge(\mathfrak{g}_{\bar 1}), \mathcal{F}_G)$ described in
$2.1$.  We have
$$
\mathcal{O}_M\subset
\widehat{\mathcal{O}}_G\stackrel{\Gamma}{\longrightarrow}
\Hom(\bigwedge(\mathfrak{g}_{\bar 1}), \mathcal{F}_G).
$$
By definition $f\in (\mathcal{O}_M)_p$ if and only if $f\in
(\widehat{\mathcal{O}}_G)_p$ and $(\mu_{G\times
H})^*(f)=\pr^*_1(f)$.
 Using (\ref{umnozh}) and (\ref{case}) we can write the last condition
 in the
following form:
\begin{equation}\label{uslovie prav invar}
f(X^r\cdot \Ad_G(g)(Y^q))(gh)=\left\{
 \begin{array}{ll}
  0, & \hbox{if}\,\,q\ne 0; \\
  f(X^r)(g), & \hbox{if}\,\, q=  0
   \end{array}
     \right.
\end{equation}
for all $g\in G$, $h\in H$, where $X^r=X_1\cdots X_r$, $X_i\in
\mathfrak{g}$, $Y^q=Y_1\cdots Y_q$, $Y_j\in \mathfrak{h}=\Lie
(H,\mathcal{O}_H)$, $r+q=p$. The supergroup $(G,\mathcal{O}_G)$ is a
split Lie supergroup, it follows by Proposition \ref{split Lie
supergroup} that $[\mathfrak{g}_{\bar 1}, \mathfrak{g}_{\bar 1}]
=0$. As we mentioned above, in this case the mapping
(\ref{isomorphism}) is an injective homomorphism of superalgebras
(not only of cosuperalgebras). Hence $s\in
\Gamma((\mathcal{O}_M)_p)$ if and only if $s\in
\Hom(\bigwedge^p(\mathfrak{g}_{\bar 1}), \mathcal{F}_G)$ and the
following condition holds:
\begin{equation}\label{uslovie prav invar rassloen}
s(X^r\wedge\Ad_G(g)( Y^q))(gh)=\left\{
 \begin{array}{ll}
  0, & \hbox{if}\,\,q\ne  0; \\
  s(X^r)(g), & \hbox{if}\,\, q= 0
   \end{array}
     \right.
\end{equation}
for all $g\in G$, $h\in H$, $X^r=X_1\wedge\cdots\wedge X_r\in
\bigwedge^r(\mathfrak{g}_{\bar 1})$, $Y^q= Y_1\wedge\cdots\wedge
Y_q\in \bigwedge^q(\mathfrak{h}_{\bar 1})$, $r+q=p$.

We may regard the sheaf $\Hom(\bigwedge(\mathfrak{g}_{\bar
1}/\mathfrak{h}_{\bar 1}), \mathcal{F}_G)$ as a sheaf of
superalgebras. The multiplication is defined as in $2.1$. Define the
injective homomorphism of sheaves of superalgebras in the following
way:
$$
\begin{array}{c}
\Upsilon: \Hom(\bigwedge(\mathfrak{g}_{\bar 1}/\mathfrak{h}_{\bar
1}), \mathcal{F}_G)\longrightarrow \Hom(\bigwedge(\mathfrak{g}_{\bar
1}), \mathcal{F}_G),\\
\Upsilon(f)(\Ad_G(g)(X))(g)=f(\overline{X})(g),
\end{array}
$$
where $X\in \bigwedge(\mathfrak{g}_{\bar 1})$, $\overline{X}$ is the
image of $X$ by the natural homomorphism
$\bigwedge(\mathfrak{g}_{\bar 1})\to \bigwedge(\mathfrak{g}_{\bar
1}/\mathfrak{h}_{\bar 1})$, $g\in G$.

Consider the composition of the injective homomorphisms of sheaves
of superalgebras
$$
\begin{array}{rl}
\bigwedge \mathcal{E}^{\psi}\stackrel{\Phi_{\wedge \psi
}}{\longrightarrow} \mathcal{F}_G\otimes
\bigwedge(\mathfrak{g}_{\bar 1}/\mathfrak{h}_{\bar 1})^*
\stackrel{\Psi_{\bigwedge(\mathfrak{g}_{\bar 1}/\mathfrak{h}_{\bar
1})^*}}{\longrightarrow} \Hom(\bigwedge(\mathfrak{g}_{\bar
1}/\mathfrak{h}_{\bar 1}),&
\mathcal{F}_G)\\
&\stackrel{\Upsilon}{\longrightarrow}
\Hom(\bigwedge(\mathfrak{g}_{\bar 1}), \mathcal{F}_G).
\end{array}
$$
Our goal now is to show that
$$
\Gamma((\mathcal{O}_M)_p)= \Upsilon\circ
\Psi_{\bigwedge(\mathfrak{g}_{\bar 1}/\mathfrak{h}_{\bar 1})^*}\circ
\Phi_{\wedge \psi }(\bigwedge^p \mathcal{E}^{\psi}),\,\,p\geq 0.
$$
This will imply our assertion.

Note that $f\in \Psi_{\bigwedge(\mathfrak{g}_{\bar
1}/\mathfrak{h}_{\bar 1})^*}\circ \Phi_{\wedge \psi }(
\bigwedge^p\mathcal{E}^{\psi})$ if and only if $f\in
\Hom(\bigwedge^p(\mathfrak{g}_{\bar 1}/\mathfrak{h}_{\bar 1}),
\mathcal{F}_G)$ and the following condition holds:
$$
f(\Ad_G(h^{-1})(\overline{X}))(gh)=f(\overline{X})(g),
$$
where $X\in \bigwedge^p(\mathfrak{g}_{\bar 1})$, $g\in G,\,\, h\in
H$. Further, $s\in \Upsilon \circ\Psi_{\bigwedge(\mathfrak{g}_{\bar
1}/\mathfrak{h}_{\bar 1})^*}\circ \Phi_{\wedge \psi }(
\bigwedge^p\mathcal{E}^{\psi})$ if and only if $s\in
\Hom(\bigwedge^p(\mathfrak{g}_{\bar 1}), \mathcal{F}_G)$ and the
following condition holds:
\begin{equation}\label{uslovie3}
s(\Ad_G(g)(X))(gh)=\left\{
                     \begin{array}{ll}
                       0, & \hbox{if $\overline{X}=\overline{0}$;} \\
                       s(\Ad_G(g)(X))(g), & \hbox{if $\overline{X}\ne \overline{0}$.}
                     \end{array}
                   \right.
\end{equation}
The conditions (\ref{uslovie prav invar rassloen}) and
(\ref{uslovie3}) are equivalent.

To complete the proof, we recall that the action
$\nu:(G,\mathcal{O}_{G})\times (M,\mathcal{O}_{M})\to
(M,\mathcal{O}_{M})$ is defined by $\nu^*(f)=\mu^*(f)$, $f\in
\mathcal{O}_{M}\subset \mathcal{O}_{G}$, and the map $\mu^*$ and the
inclusion $\mathcal{O}_{G\times M}\hookrightarrow
\mathcal{O}_{G\times G}$ preserve the chosen
$\mathbb{Z}$-gradings.$\Box$

\medskip

Let us formulate the general result concerning a complex homogeneous
split supermanifold.

\medskip

\t\label{split} {\it Let $(G,\mathcal{O}_G)$ be a complex Lie
supergroup with the Lie superalgebra $\mathfrak{g}=
\mathfrak{g}_{\bar 0} \oplus \mathfrak{g}_{\bar 1}$. If
$[\mathfrak{g}_{\bar 1}, \mathfrak{g}_{\bar 1}] =0$ then all
$(G,\mathcal{O}_G)$-homogeneous supermanifolds $(M,\mathcal{O}_M)$
are split supermanifolds. Moreover, the sheaf $\mathcal{O}_M$ is
isomorphic to $\bigwedge \mathcal{E}^{\psi}$, where
$\mathcal{E}^{\psi}$ is the sheaf of sections of the homogeneous
vector bundle $\mathbf{E}^{\psi}$, which corresponds to the
representation $\psi: H \to GL((\mathfrak{g}_{\bar
1}/\mathfrak{h}_{\bar 1})^*)$ given by
$$
\psi(h)(v)(X+ \mathfrak{h}_{\bar 1}):= v(\Ad_G(h^{-1})(X)+
\mathfrak{h}_{\bar 1}),\,\,\text{for}\,\,h\in H, \,X\in
\mathfrak{g}_{\bar 1},\,\, v\in (\mathfrak{g}_{\bar
1}/\mathfrak{h}_{\bar 1})^*.
$$

 Conversely, if a complex homogeneous supermanifold
$(M,\mathcal{O}_M)$ is split then there is a Lie supergroup
$(G,\mathcal{O}_G)$ with $[\mathfrak{g}_{\bar 1}, \mathfrak{g}_{\bar
1}] =0$, where $\mathfrak{g} = \mathfrak{g}_{\bar 0}\oplus
\mathfrak{g}_{\bar 1} = Lie (G,\mathcal{O}_G)$, such that
$(G,\mathcal{O}_G)$ acts on $(M,\mathcal{O}_M)$ transitively.

}

\medskip

\noindent{\it Proof.} The theorem follows from Propositions
\ref{split Lie supergroup}, \ref{quptient of split Lie supergroups}
and Lemma \ref{split action}.$\Box$

\medskip

Let us prove for example that the complex projective superspace
$\mathbb{CP}^{1|2}$ is split. It is isomorphic to the coset space
$\GL_{2|1}(\mathbb{C})/(P,\mathcal{O}_P)$, where
$$
\GL_{2|1}(\mathbb{C})=\left(
                        \begin{array}{ccc}
                          * & * & \checkmark \\
                          * & * & \checkmark \\
                          \checkmark & \checkmark & * \\
                        \end{array}
                      \right),\,\,\,
(P,\mathcal{O}_P)= \left(
                        \begin{array}{ccc}
                          * & * & \checkmark \\
                          0 & * & \checkmark \\
                          0 & \checkmark & * \\
                        \end{array}
                      \right).
$$
Here $*$ are even coordinates and $\checkmark$ are odd coordinates.
It is easy to see that
$$
\GL_{2|1}(\mathbb{C})/(P,\mathcal{O}_P) \simeq
(G',\mathcal{O}_{G'})/(P',\mathcal{O}_{P'}),
$$
 where
$$
(G',\mathcal{O}_{G'})=\left(
                        \begin{array}{ccc}
                          * & * & 0 \\
                          * & * & 0 \\
                          \checkmark & \checkmark & * \\
                        \end{array}
                      \right),\,\,\,
(P',\mathcal{O}_{P'})= \left(
                        \begin{array}{ccc}
                          * & * & 0 \\
                          0 & * & 0 \\
                          0 & \checkmark & * \\
                        \end{array}
                      \right).
$$
Let $\mathfrak{g}'=\operatorname{Lie}(G',\mathcal{O}_{G'})$. Then
$$
\mathfrak{g}'_{\bar 1}= \left\{\left(
                        \begin{array}{ccc}
                          0 & 0 & 0 \\
                          0 & 0 & 0 \\
                          a & b & 0 \\
                        \end{array}
                      \right),\quad  a,b\in \mathbb{C} \right\}.
$$
We see that $[\mathfrak{g}'_{\bar 1}, \mathfrak{g}'_{\bar
1}]=\{0\}$. By Theorem \ref{split} we get that $\mathbb{CP}^{1|2}$
is split.

We close this section by mentioning some results about non-split
supermanifolds. The first example of a non-split supermanifold was
published in \cite{Green}; this is the quadric in the projective
superplane $\mathbb{CP}^{2|2}$. In \cite{man} four series of
supermanifolds of flags were constructed corresponding to four
series of classical linear Lie superalgebras. In \cite{LeBrun} it
was proved that all split complex supermanifolds whose reduction is
projective algebraic are projective (that is, embeddable in a
complex projective superspace.) Penkov and Skornyakov \cite{Penkov}
found necessary and sufficient conditions for a supermanifold of
flags to be projective. More precisely, they showd that almost all
such supermanifolds are not projective. From these two results it
follows that supermanifolds of flags are mostly non-split.

In \cite{onipisp} it was proved that the isotropic
super-Grassmannian of maximal type
$\mathbb{IG}\!\operatorname{r}_{n|n,s|t}(\mathbb{C})$ associated
with an odd bilinear form is non-split whenever $t\geq 1$ and $s\geq
2$. In \cite{oniosp} the complete solution of the problem was given
for the isotropic super-Grassmannian of maximal type associated with
an even bilinear form. Note that the method of \cite{oniosp} and
\cite{onipisp} can be used for all series of flag supermanifolds.

In \cite{onigrassmann} the problem of classifying all homogeneous
complex supermanifolds whose reduction is the complex Grassmannian
$\mathrm{Gr}_{n|k}$ was studied. Under the assumption that the odd
isotropy representation is irreducible and under certain
restrictions on $(n|k)$, it was proved that the only non-split
supermanifold of this sort is the $\Pi$-symmetric super-Grassmannian
constructed by Manin \cite{man}.

The problem of classification of non-split supermanifolds having as
retract the split supermanifold $(M,\Omega)$, where $\Omega$ is the
sheaf of holomorphic forms on a given complex manifold $M$ of
dimension $> 1$, was studied in \cite{oniCOT}. In the case when $M$
is an irreducible compact Hermitian symmetric space, the complete
classification of non-split supermanifolds with retract $(M,\Omega)$
was given.

\smallskip

\noindent  {\bf Acknowledgment.} The author is grateful to A.L.
Onishchik, P. Heinzner, A. T. Huckleberry and anonymous referees for
useful comments.

\smallskip

\noindent {\textsc{Ruhr-Universit\"{a}t Bochum,
Universit\"{a}tsstra{\ss}e 150, 44780 Bochum, Germany;}}

 \noindent {\textsc{Tver State University, Zhelyabova 33, 170 000
Tver, Russia}};

 \noindent {\emph{E-mail address:}
\verb"VishnyakovaE@googlemail.com"}


\begin{thebibliography}{99}

\bibitem{BagSta} {\it Baguis P.,  Stavracou T.} Normal Lie subsupergroups and non-abelian
supercircles. International Journal of Mathematics and Mathematical
Sciences Volume 30 (2002), Issue 10, Pages 581-591.


\bibitem{BL} {\it Berezin F.A., Leites D.A.} Supermanifolds. Soviet
Math. Dokl. 16, 1975, 1218-1222.

\bibitem{Bern} {\it Deligne P., Morgan J.W.} Notes on supersymmetry
(following Joseph Bernstein), Quantum Fields and Strings: A Course
for Mathematicians, Vols. 1,2 (Princeton, NJ, 1996/1997), 41-97.
American Mathematical Society. Providence, R.I. 1999.

\bibitem{Fio_Varad} {\it   Fioresi R., Lledo M.A., Varadarajan V. S.} The
super Minkowski and conformal space times, JMP, 48, no. 11, pg.
113505, 2007.

\bibitem{Green} {\it   Green P.} On holomorphic graded manifolds,  Proc. Amer. Math. Soc.  85  (1982), no. 4, 587--590.


\bibitem{Kostant} {\it Kostant B.} Graded manifolds, graded Lie
theory, and prequantization. Lecture Notes in Mathematics 570.
Berlin e.a.: Springer-Verlag,  1977. P. 177-306.

\bibitem{kosz} {\it  Koszul J.L.} Graded manifolds and graded Lie
algebras. Proceeding of the International Meeting on Geometry and
Physics (Bologna), Pitagora, 1982, pp 71-84.

\bibitem{LeBrun} {\it LeBrun C., Poon Y. S., Wells R. O.,
Jr.} Projective embeddings of complex supermanifolds.  Comm. Math.
Phys.  126  (1990),  no. 3, 433--452.


\bibitem{ley} {\it Leites D.A.} Introduction to the theory of supermanifolds.
Russian Math. Surveys 35 (1980), 1-64.

\bibitem{man} {\it Manin Yu. I.} Gauge field theory and complex
geometry. Grundlehren der Mathematischen Wissenschaften, 289,
Springer-Verlag, Berlin, 1997.

\bibitem{Mol} {\it Molotkov V.} Infinite-dimensional $\mathbb{Z}^k_2$-supermanifolds.
International center for theoretical physics, Preprint ¹ IC/84/183,
1984, http://ccdb4fs.kek.jp/cgi-bin/img/allpdf?198506284.

\bibitem{onipi} {\it Onishchik A.L.} Flag supermanifolds, their
automorphisms und deformations. The Sophus Lie Memorial conference
(Oslo, 1992), 289-302, Scand. Univ. Press, Oslo, 1994.

\bibitem{onigrassmann} {\it Onishchik A.L.} Homogeneous supermanifolds over
Grassmannians,  J. Algebra  313  (2007),  no. 1, P. 320--342.

\bibitem{oniCOT} {\it Onishchik A.L.} Non-split supermanifolds associated with the
cotangent bundle. Universit\'{e} de Poitiers, D\'{e}partement de
Math., N 109. Poitiers, 1997.

\bibitem{oniosp} {\it Onishchik A.L., Serov A.A.} Vector fields and deformations of
isotropic super-Grassmannians of maximal type. Lie Groups and Lie
Algebras: E.B. Dynkin's Seminar, AMS Transl. Ser. 2. V. 169.
Providence: AMS,  1995. P. 75-90.

\bibitem{onipisp} {\it Onishchik A.L., Serov A.A.} On isotropic super-Grassmannians
of maximal type associated with an odd bilinear form. E.
Schr\"{o}dinger Inst. for Math. Physics, Preprint No. 340. Vienna,
1996.

\bibitem{Penkov}{\it Penkov I.B., Skornyakov I.A. } Projectivity and $\mathrm{D}$-affineness of flag
supermanifolds, Russ. Math. Surv. 40, (1987) P. 233-234.

\bibitem{Scheu} {\it Scheunert M.}  The theory of Lie
superalgebras. Lectures Notes in Mathematics 716, Springer, Berlin
1979.

\bibitem{TSH} {\it Tsalenko M.S., Shulgeifer E.G.}  Foundations of the
theory of categories. Nauka, Moscow 1974 (in Russian).

\bibitem{Var} {\it Varadarajan V. S.} Supersymmetry for
mathematicians: an introduction, AMS, Courant lecture notes, Vol.
11, 2004.

\bibitem{V} {\it Vishnyakova E. G.}  On the structure of complex homogeneous supermanifolds. arXiv:0811.2581 (November 2008).

\bibitem{Westra} {\it Westra D.B.}  Superrings and supergroups. PhD
thesis, Universit\"{a}t Wien, 2009.



\end{thebibliography}
\end{document}